

\magnification\magstep1
\openup1\jot

\def\bquote{``\thinspace\thinspace}
\def\equote{\thinspace\thinspace''}

\hsize=11.5cm    
\vsize=18cm      
\parindent=12pt   \parskip=5pt    

\hoffset=.5cm  
\voffset=.8cm  

\pretolerance=500 \tolerance=1000  \brokenpenalty=5000

\catcode`\@=11

\font\eightrm=cmr8         \font\eighti=cmmi8
\font\eightsy=cmsy8        \font\eightbf=cmbx8
\font\eighttt=cmtt8        \font\eightit=cmti8
\font\eightsl=cmsl8        \font\sixrm=cmr6
\font\sixi=cmmi6           \font\sixsy=cmsy6
\font\sixbf=cmbx6

\font\tengoth=eufm10
\font\eightgoth=eufm8  
\font\sevengoth=eufm7      
\font\sixgoth=eufm6        \font\fivegoth=eufm5

\skewchar\eighti='177 \skewchar\sixi='177
\skewchar\eightsy='60 \skewchar\sixsy='60

\newfam\gothfam           \newfam\bboardfam

\def\tenpoint{
  \textfont0=\tenrm \scriptfont0=\sevenrm \scriptscriptfont0=\fiverm
  \def\rm{\fam\z@\tenrm} \textfont1=\teni \scriptfont1=\seveni
  \scriptscriptfont1=\fivei
  \def\oldstyle{\fam\@ne\teni}\let\old=\oldstyle \textfont2=\tensy
  \scriptfont2=\sevensy \scriptscriptfont2=\fivesy
  \textfont\gothfam=\tengoth \scriptfont\gothfam=\sevengoth
  \scriptscriptfont\gothfam=\fivegoth \def\goth{\fam\gothfam\tengoth}
  
  \textfont\itfam=\tenit
  \def\it{\fam\itfam\tenit}
  \textfont\slfam=\tensl
  \def\sl{\fam\slfam\tensl}
  \textfont\bffam=\tenbf \scriptfont\bffam=\sevenbf
  \scriptscriptfont\bffam=\fivebf
  \def\bf{\fam\bffam\tenbf}
  \textfont\ttfam=\tentt
  \def\tt{\fam\ttfam\tentt}
  \abovedisplayskip=12pt plus 3pt minus 9pt
  \belowdisplayskip=\abovedisplayskip
  \abovedisplayshortskip=0pt plus 3pt
  \belowdisplayshortskip=4pt plus 3pt
  \smallskipamount=3pt plus 1pt minus 1pt
  \medskipamount=6pt plus 2pt minus 2pt
  \bigskipamount=12pt plus 4pt minus 4pt
  \normalbaselineskip=12pt
  \setbox\strutbox=\hbox{\vrule height8.5pt depth3.5pt width0pt}
  \let\bigf@nt=\tenrm       \let\smallf@nt=\sevenrm
  \normalbaselines\rm}

\def\eightpoint{
  \textfont0=\eightrm \scriptfont0=\sixrm \scriptscriptfont0=\fiverm
  \def\rm{\fam\z@\eightrm}
  \textfont1=\eighti  \scriptfont1=\sixi  \scriptscriptfont1=\fivei
  \def\oldstyle{\fam\@ne\eighti}\let\old=\oldstyle
  \textfont2=\eightsy \scriptfont2=\sixsy \scriptscriptfont2=\fivesy
  \textfont\gothfam=\eightgoth \scriptfont\gothfam=\sixgoth
  \scriptscriptfont\gothfam=\fivegoth
  \def\goth{\fam\gothfam\eightgoth}
  
  \textfont\itfam=\eightit
  \def\it{\fam\itfam\eightit}
  \textfont\slfam=\eightsl
  \def\sl{\fam\slfam\eightsl}
  \textfont\bffam=\eightbf \scriptfont\bffam=\sixbf
  \scriptscriptfont\bffam=\fivebf
  \def\bf{\fam\bffam\eightbf}
  \textfont\ttfam=\eighttt
  \def\tt{\fam\ttfam\eighttt}
  \abovedisplayskip=9pt plus 3pt minus 9pt
  \belowdisplayskip=\abovedisplayskip
  \abovedisplayshortskip=0pt plus 3pt
  \belowdisplayshortskip=3pt plus 3pt
  \smallskipamount=2pt plus 1pt minus 1pt
  \medskipamount=4pt plus 2pt minus 1pt
  \bigskipamount=9pt plus 3pt minus 3pt
  \normalbaselineskip=9pt
  \setbox\strutbox=\hbox{\vrule height7pt depth2pt width0pt}
  \let\bigf@nt=\eightrm     \let\smallf@nt=\sixrm
  \normalbaselines\rm}

\tenpoint

\def\pc#1{\bigf@nt#1\smallf@nt}         \def\pd#1 {{\pc#1} }

\catcode`\;=\active
\def;{\relax\ifhmode\ifdim\lastskip>\z@\unskip\fi
\kern\fontdimen2  -1.2 \fontdimen3 \string;}

\catcode`\:=\active
\def:{\relax\ifhmode\ifdim\lastskip>\z@\unskip\fi\penalty\@M\ \fi\string:}

\catcode`\!=\active
\def!{\relax\ifhmode\ifdim\lastskip>\z@
\unskip\fi\kern\fontdimen2  -1.1 \fontdimen3 \string!}

\catcode`\?=\active
\def?{\relax\ifhmode\ifdim\lastskip>\z@
\unskip\fi\kern\fontdimen2  -1.1 \fontdimen3 \string?}

\catcode`\«=\active
\def«{\raise.4ex\hbox{%
 $\scriptscriptstyle\langle\!\langle$}}

\catcode`\»=\active
\def»{\raise.4ex\hbox{%
 $\scriptscriptstyle\rangle\!\rangle$}}

\frenchspacing

\def\raggedbottom{\topskip 10pt plus 36pt\r@ggedbottomtrue}

\def\pointir{\unskip . --- \ignorespaces}

\def\Medbreak{\vskip-\lastskip\medbreak}

\long\def\th#1 #2\enonce#3\endth{
   \Medbreak\noindent
   {\pc#1} {#2\unskip}\pointir{\it #3}\smallskip}

\def\proof{\vskip-\lastskip\smallskip\noindent
 {\it Proof} : }

\def\decale#1{\smallbreak\hskip 28pt\llap{#1}\kern 5pt}
\def\decaledecale#1{\smallbreak\hskip 34pt\llap{#1}\kern 5pt}
\def\puce{\smallbreak\hskip 6pt{$\scriptstyle\bullet$}\kern 5pt}

\def\eqalign#1{\null\,\vcenter{\openup\jot\m@th\ialign{
\strut\hfil$\displaystyle{##}$&$\displaystyle{{}##}$\hfil
&&\quad\strut\hfil$\displaystyle{##}$&$\displaystyle{{}##}$\hfil
\crcr#1\crcr}}\,}

\catcode`\@=12

\showboxbreadth=-1  \showboxdepth=-1

\mathcode`A="7041 \mathcode`B="7042 \mathcode`C="7043 \mathcode`D="7044
\mathcode`E="7045 \mathcode`F="7046 \mathcode`G="7047 \mathcode`H="7048
\mathcode`I="7049 \mathcode`J="704A \mathcode`K="704B \mathcode`L="704C
\mathcode`M="704D \mathcode`N="704E \mathcode`O="704F \mathcode`P="7050
\mathcode`Q="7051 \mathcode`R="7052 \mathcode`S="7053 \mathcode`T="7054
\mathcode`U="7055 \mathcode`V="7056 \mathcode`W="7057 \mathcode`X="7058
\mathcode`Y="7059 \mathcode`Z="705A

\font\ss=cmss10
\font\bigss=cmss11

\def\qp{{\bf Q}_p}

\def\br#1{\hbox{{\ss Br}}(#1)}

\def\ache{\hbox{\ss H}}

\def\kbar{{\overline k}}
\def\Fbar{{\overline F}}

\def\kbaretoile{{\kbar{}^{\times}}}

\def\Finfty{{F_{\!\infty}}}
\def\Finftyetoile{{F_{\!\infty}^\times}}

\def\det{\mathop{\hbox{d{\'e}t}}\nolimits}
\def\Gal{\mathop{\hbox{\ss Gal}}\nolimits}

\def\br(#1){\hbox{{\ss B}}(#1)}

\def\Hom{\mathop{\hbox{\ss Hom}}\nolimits}

\def\End{\mathop{\rm End}\nolimits}
\def\Card{\mathop{\rm Card}\nolimits}
\def\Spec{\mathop{\hbox{\ss Spec}}\nolimits}

\def\F{\mathord{\bf F}}
\def\P{\mathord{\bf P}}
\def\M{\mathord{\bf M}}
\def\Z{\mathord{\bf Z}}
\def\R{\mathord{\bf R}}
\def\Q{\mathord{\bf Q}}

\def\ogoth{{\goth o}}
\def\pgoth{{\goth p}}

\def\hfl#1#2#3{\smash{\mathop{\hbox to#3{\rightarrowfill}}\limits
^{\textstyle#1}_{\textstyle#2}}}
\def\gfl#1#2#3{\smash{\mathop{\hbox to#3{\leftarrowfill}}\limits
^{\textstyle#1}_{\textstyle#2}}}

\def\qed{\raise -2pt\hbox{\vrule\vbox to 10pt{\hrule width 4pt
                 \vfill\hrule}\vrule}}
\def\cqfd{\unskip\penalty 500\quad\qed\medbreak}

\def\phi{\varphi}

\def\Ker{\mathop{\rm Ker}\nolimits}
\def\Coker{\mathop{\rm Coker}\nolimits}
\def\Id{\mathop{\rm Id}\nolimits}

\def\droite#1{\,\hfl{#1}{}{8mm}\,}

\def\diagram#1{\def\normalbaselines{\baselineskip=0pt\lineskip=5pt}
\matrix{#1}}
\def\vfl#1#2#3{\llap{$\textstyle #1$}
\left\downarrow\vbox to#3{}\right.\rlap{$\textstyle #2$}}

\newcount\refno

\long\def\ref#1:#2<#3>{                                        
\global\advance\refno by1\par\noindent                              
\llap{[{\bf\number\refno}]\ }{#1}\pointir{\it #2} #3\goodbreak }

\def\citer#1(#2){[{\bf\number#1}\if#2\empty\relax\else,\ #2\fi]}

\def\fleche{\rightarrow}

\newcount\numerodesection \numerodesection=1
\def\section#1{\bigbreak
 {\bf\number\numerodesection.\ \ #1}\nobreak\medskip
 \advance\numerodesection by1}

\newcount\numeroderemarque
\def\remarque{\advance\numeroderemarque by1\smallbreak
{\it Remarque\/}\ \number\numeroderemarque~:}

\newcount\formuleno
\def\numeroter{\global\advance\formuleno by1
 \leqno{(\oldstyle\number\formuleno)}}
\def\formule#1{$(\oldstyle\number#1)$}

\def\div{\mathop{\rm div}\nolimits}

\def\pr{\mathop{\rm pr}\nolimits}
\def\somme{\mathop{\smash{\raise 2pt\hbox{$\sum$}}}\limits}
\def\sommedir{\mathop{\oplus}\limits}

\def\Tor{\mathop{\rm Tor}\nolimits}
\def\zero{\{0\}}
\def\one{\{1\}}
\def\long{\mathop{\rm long}\nolimits}

\def\Arg{\mathop{\rm Arg}\nolimits}

\def\symbol(#1,#2){\left\{{#1},{#2}\right\}}
\def\qfclass(#1){\left\langle{#1}\right\rangle}
\def\commut(#1,#2){\lgroup#1,#2\rgroup}
\def\kay{\hbox{\ss K}}
\def\vv{\hbox{\ss V}}

\def\anought{\hbox{\ss A}_0}
\def\kaytwo{\hbox{\ss K}_2}
\def\kayenn{\hbox{\ss K}_n}
\def\derivedsg{\hbox{\ss D}}
\def\witt{\hbox{\ss W}}
\def\witthat{\hbox{$\hat{\hbox{\ss W}}$}}
\def\quadsp{\hbox{\ss Q}}
\def\Fstar{F^\times}
\def\Lstar{L^\times}
\def\Fprime{F^\prime}
\def\Aprime{A^\prime}
\def\kstar{k^\times}
\def\Fpstar{{\bf F}_{\!p}^\times}
\def\Fsep{{\overline F}}
\def\Fsepstar{{\overline F}{}^\times}
\def\Fab{{\tilde F}}
\def\Qp{{\bf Q}_p}
\def\Fp{{\bf F}_p}
\def\Qbar{\bar{\Q}}
\def\C{{\bf C}}
\def\N{{\bf N}}
\def\Qpstar{{\bf Q}_p^\times}
\def\Zpstar{{\bf Z}_p^\times}
\def\inv{^{-1}}
\def\xinv{x\inv}
\def\SL{\mathop{\bf SL}\nolimits}
\def\Sp{\mathop{\bf Sp}\nolimits}
\def\GL{\mathop{\bf GL}\nolimits}
\def\diag{\mathop{\rm diag}\nolimits}
\def\pr{\mathop{\rm pr}\nolimits}
\def\dlog{\mathop{\rm dlog}\nolimits}
\def\derlog(#1){{d#1\over#1}}
\def\normressym(#1,#2)_#3{\displaystyle\left({#1,#2\over#3}\right)}
\def\cupprod{\mathrel{\lower3pt\hbox{$\smile$}}}

\newbox\bibbox
\setbox\bibbox\vbox{\bigskip
\centerline{---$*$---$*$---}
\bigbreak
\centerline{{\pc BIBLIOGRAPHICAL REFERENCES}}
We content ourselves with a very short bibliography.  Many of the
papers are available at {\tt arxiv.org}, or {\tt gallica.bnf.fr}, or
{\tt gdz.sub.uni-goettingen.de}, or {\tt www.numdam.org}.


\def\mr#1.{{\sevenrm MR#1.}}

\ref {\pc AKHTAR} (Reza) :
Milnor $K$-theory of smooth varieties,
<$K$-Theory {32} (2004), no.~3, 269--291.
\mr2114169.>
\newcount\akhtar \global\akhtar=\refno

\ref {\pc BASS} (Hyman) :
$K_2$ des corps globaux,
<[d'apr{\`e}s J. Tate, H. Garland]
S{\'e}minaire Bourbaki, 23\raise3pt\hbox{\eightrm{\`e}me}
ann{\'e}e (1970/1971), Exp.\ 394, pp.\
233--255. Lecture Notes in Math., Vol.\ 244, Springer, Berlin,
1971. \mr0422211.>
\newcount\bass \global\bass=\refno

\ref {\pc BASS} (Hyman) and {\pc TATE} (John) :
The Milnor ring of a global field,
<in Algebraic $K$-theory,
II~: ``\thinspace Classical\thinspace'' algebraic $K$-theory and
connections with arithmetic
(Proc.\ Conf., Battelle Memorial Inst., Seattle, Wash., 1972),
pp. 489--501. Lecture Notes in Math., Vol. 342, Springer, Berlin,
1973. \mr0442061.>
\newcount\basstate \global\basstate=\refno

\ref {\pc BEILINSON} (Alexander) :
Higher regulators and values of $L$-functions.
<Current problems in mathematics, Vol. 24, 181--238,
Itogi Nauki i Tekhniki,
Akad.\  Nauk SSSR, Vsesoyuz.\ Inst.\ Nauchn.\ i Tekhn.\ Inform.,
Moscow, 1984. \mr0760999.> 
\newcount\beilinson \global\beilinson=\refno

\ref {\pc BLOCH} (Spencer) :
Algebraic $K$-theory and classfield theory for arithmetic surfaces,
<Ann.\ of Math.\ (2) 114 (1981), no.~2, 229--265. \mr0632840.>
\newcount\blochcft \global\blochcft=\refno

\ref {\pc BLOCH} (Spencer) :
Higher regulators, algebraic $K$-theory, and zeta functions of elliptic curves.
<CRM Monograph Series, 11.
American Mathematical Society, Providence, RI, 2000. \mr1760901.>
\newcount\bloch \global\bloch=\refno

\ref {\pc BLOCH} (Spencer) \& {\pc GRAYSON} (Daniel) :
$K\sb 2$ and $L$-functions of elliptic curves: computer calculations.
<applications of algebraic $K$-theory to algebraic geometry and number theory, Part I, II (Boulder, Colo., 1983), 79--88,
Contemp.\ Math., 55,
Amer.\ Math.\ Soc., Providence, RI, 1986. \mr0862631.>
\newcount\blochgrayson \global\blochgrayson=\refno


\ref {\pc BLOCH} (Spencer) \& {\pc  KATO} (Kazuya) :
$p$-adic {\'e}tale cohomology,
<Inst.\ Hautes {\'E}tudes sci.\ Publ.\ Math.\ No.\ 63 (1986), 107--152.
\mr0849653.>
\newcount\bgk \global\bgk=\refno

\ref {\pc CHASE} (Stephen) \& {\pc WATERHOUSE} (William) :
Moore's theorem on uniqueness of reciprocity laws.
<Invent.\ math. 16 (1972), 267--270. \mr0311623.>
\newcount\chasewater \global\chasewater=\refno

\ref {\pc COLLIOT-\pc TH{\'E}L{\`E}NE} (Jean-Louis) :
Cohomologie galoisienne des corps valu{\'e}s discrets henseliens,
d'apr{\`e}s K. Kato et S. Bloch,
<in Algebraic $K$-theory and its applications
(Trieste, 1997), 120--163, World Sci. Publishing, River Edge, NJ,
1999.  {\tt www.math.u-psud.fr/\char`~colliot/liste-publi.html},
\mr1715874.>
\newcount\colliot \global\colliot=\refno

\ref{\pc DIEUDONN{\'E}} (Jean) :
A. Grothendieck's early work (1950--1960).
<$K$-Theory 3 (1989), no.~4, 299--306.  \mr1047190.>
\newcount\dieudonne \global\dieudonne=\refno

\ref {\pc DOKCHITSER} (Tim), {\pc {}DE \pc JEU} (Rob) \& {\pc ZAGIER}
(Don) : 
Numerical verification of Beilinson's conjecture for\/ $K_2$ of
hyperelliptic curves. 
<to appear in Compositio math.\ Preprint math.AG/0405040>
\newcount\ddjz \global\ddjz=\refno

\ref {\pc KAHN} (Bruno) :
La conjecture de Milnor,
<Ast\'erisque No. 245 (1997), Exp.\ No.\ 834, 5, 379--418.
{\tt www.math.uiuc.edu/K-theory/0210}, \mr1627119.>
\newcount\kahn \global\kahn=\refno

\ref {\pc KATO} (Kazuya) :
Symmetric bilinear forms, quadratic forms
and Milnor $K$-theory in characteristic two,
<Invent.\ math.\ 66 (1982),  no. 3, 493--510.  \mr0662605.>
\newcount\kato \global\kato=\refno

\ref {\pc KAROUBI} (Max) :
Rapport sur la $K$-th{\'e}orie (1956--1997). 
<Development of mathematics 1950--2000,
635--653, Birkh{\"a}user, Basel, 2000. \mr1796854.>
\newcount\karoubi \global\karoubi=\refno

\ref  {\pc LICHTENBAUM} (Stephen) :
Values of zeta-functions, {\'e}tale cohomology, and algebraic
$K$-theory,
<in Algebraic $K$-theory,
II~: ``\thinspace Classical\thinspace'' algebraic $K$-theory and
connections with arithmetic
(Proc.\ Conf., Battelle Memorial Inst., Seattle, Wash., 1972),
pp. 489--501. Lecture Notes in Math., Vol. 342, Springer, Berlin,
1973.  \mr0406981.>
\newcount\licht \global\licht=\refno

\ref {\pc MATSUMOTO} (Hideya) :
Sur les sous-groupes arithm{\'e}tiques des groupes semi-simples
d{\'e}ploy{\'e}s,
<Ann.\ Sci.\ {\'E}cole Norm.\ Sup.\ (4) 2 (1969), 1--62.  \mr0240214.>
\newcount\matsumoto \global\matsumoto=\refno

\ref {\pc MAZUR} (Barry) \& {\pc WILES} (Andrew) :
Class fields of abelian extensions of\/ ${\bf Q}$,
<Invent.\ math.\ 76 (1984), no.\ 2, 179--330.  \mr0742853.>
\newcount\mazurwiles \global\mazurwiles=\refno

\ref {\pc MILNOR} (John) :
Algebraic $K$-theory and quadratic forms,
<Invent.\ math.\ 9 (1969/1970) 318--344.  
\mr0260844.>
\newcount\milnor \global\milnor=\refno

\ref {\pc MILNOR} (John) :
Introduction to algebraic $K$-theory.
<Annals of Mathematics Studies, No.\ 72. Princeton
University Press, Princeton, N.J., 1971. xiii+184 pp.  
\mr0349811.>
\newcount\mbook \global\mbook=\refno

\ref {\pc MOORE} (Calvin) :
Group extensions of $p$-adic and adelic linear groups,
<Inst.\ Hautes {\'E}tudes sci.\ Publ.\ math.\ No.\ 35,
(1968) 157--222. \mr0244258.>
\newcount\moore \global\moore=\refno

\ref {\pc MOREL} (Fabien) :
Voevodsky's proof of Milnor's conjecture,
<Bull.\ Amer.\ Math.\ Soc.\ (N.S.) 35 (1998), no.\ 2, 123--143.  
{\tt www.ams.org/bulletin}, \mr1600334.>
\newcount\morel \global\morel=\refno

\ref {\pc RASKIND} (Wayne) :
 Abelian class field theory of arithmetic schemes,
<in $K$-theory and algebraic geometry: connections
with quadratic forms and division algebras (Santa Barbara, CA, 1992),
85--187, Proc.\ Sympos.\ Pure Math., 58, Part 1, Amer.\ Math.\ Soc.,
Providence, RI, 1995.  \mr1327282.>
\newcount\raskind \global\raskind=\refno

\ref {\pc RUBIN} (Karl) :
Appendix,
<in Serge Lang, {\it Cyclotomic fields\/} I {\it and\/} II.
Combined second edition.  Graduate Texts in Mathematics, 121.
Springer-Verlag, New York, 1990. xviii+433 pp.
\mr1029028.>
\newcount\rubin \global\rubin=\refno

\ref {\pc SAITO} (Shuji) :
Class field theory for curves over local fields.
<J. Number Theory 21 (1985), no.~1, 44--80. \mr0804915.>
\newcount\saito \global\saito=\refno

\ref {\pc SOMEKAWA} (Matsuro) :
On Milnor $K$-groups attached to semi-abelian varieties,
<$K$-Theory 4 (1990), no.~2, 105--119.
\mr1081654.>
\newcount\somekawa \global\somekawa=\refno

\ref {\pc SOUL{\'E}} (Christophe) :
$K_2$ et le groupe de Brauer (d'apr{\`e}s A. S. Merkurjev et
A. A. Suslin),
<Ast{\'e}risque, 105-106, Soc.~Math.~France,
Paris, 1983.  \mr0728982.>
\newcount\soule \global\soule=\refno

\ref {\pc STEINBERG} (Robert) :
G{\'e}n{\'e}rateurs,
relations et rev{\^e}tements de groupes alg{\'e}briques,
<in Colloq.\ Th{\'e}orie des Groupes alg{\'e}briques (Bruxelles, 1962)
Librairie universitaire, Louvain; Gauthier-Villars,
Paris.  \mr0153677.>
\newcount\steinberg \global\steinberg=\refno

\ref {\pc TATE} (John) :
Symbols in arithmetic,
<in Actes du Congr{\`e}s international des math{\'e}maticiens (Nice,
1970), tome 1, pp.\ 201--211. Gauthier-Villars, Paris, 1971.
\mr0422212.>
\newcount\tatesym \global\tatesym=\refno

\ref {\pc TATE} (John) :
Relations between $K\sb{2}$ and Galois cohomology,
<Invent.\ math.\ 36 (1976), 257--274.  \mr0429837.>
\newcount\tatekii \global\tatekii=\refno


\ref{\pc WEIBEL} (Charles) :
The Development of Algebraic K-theory before 1980.
<{\tt http://www.math.rutgers.edu/\char`~weibel/}.>
\newcount\weibel \global\weibel=\refno

\ref {\pc YOSHIDA} (Teruyoshi) :
Finiteness theorems in the class field theory of varieties over local
fields,
<J.\ Number Theory 101 (2003), no.~1, 138--150. \mr1979656.>
\newcount\yoshida \global\yoshida=\refno
}

\centerline{\bf Some aspects of the functor\/ $K_2$ of fields}
\medskip
\centerline{Chandan Singh Dalawat}
\bigskip

We review connections between the group $\kaytwo$ of a field and universal
central extensions, quadratic forms, central simple algebras, differential
forms, abelian extensions, abelian coverings, explicit reciprocity laws,
special values of zeta functions, and special values of L-functions. No
proofs, minimal bibliography.

%

\section{The birth of $\kaytwo$ \citer\dieudonne(), \citer\karoubi(),
\citer\weibel().} 

It was in 1956 that A. Grothendieck had the idea of assigning to every scheme
$X$ a group $\kay(X)$; the letter $\kay$ stands for {\it Klassen\/}, and
indeed the elements of $\kay(X)$ are equivalence classes of certain objects
associated to $X$.  The need for such a group arose in his generalisation of
the Riemann-Roch theorem of F.~Hirzebruch from complex analytic geometry to a
relative version in abstract algebraic geometry.

Let us recall the defintion of this group when $X=\Spec(A)$ is an affine
scheme, i.e.~the spectrum of a commutative ring $A$.  Isomorphism classes of
finitely generated projective $A$-modules form a commutative monoid $M$ whose
law of composition is given by taking the the direct sum.  There arise in
practice a number of morphisms from M into commutative groups; Grothendieck
defined $\kay(X)$ --- also denoted $\kay(A)$ in this affine case --- as the
universal object amongst them.  Thus, $\kay(A)$ is simply the group of
differences of the monoid $M$, just as $\Z$ is the group of difference of the
monoid $\N$; it is a powerful invariant of the ring $A$.  For example, when
$A$ is the ring of integers in a finite extension of $\Q$, the torsion
subgroup of $\kay(A)$ is the group of ideal classes of $A$.  When $A$ is a
field $F$, the rank or dimension gives an isomorphism $K(F)\rightarrow\Z$.  


H. Bass realised that a pre-Grothendieck construction of J. H. C. Whitehead
--- the quotient of $\GL(A)=\lim\GL_n(A)$ by its derived subgroup ---
should be the right definition of the $\kay_1$ of a ring $A$.  When $A$ is a
field $F$, the determinant gives an isomorphism $\kay_1(F)\rightarrow
F^\times$. 

The next step was taken by J. Milnor who defined the group $\kaytwo(A)$ and a
cup-product map $\kay_1(A)\times\kay_1(A)\rightarrow\kaytwo(A)$ for every ring
$A$; it turned out to be the homology group $\ache_2(\derivedsg(\GL(A)),\Z)$,
where $\derivedsg(\phantom{G})$ is the derived subgroup.  A deep theorem of
H. Matsumoto gave a presentation of the $\kaytwo$ of a field $F$.  We have
chosen to take this presentation as our starting point to recount a few places
where the group $\kaytwo$ makes its appearance.

We have had to leave out more than we have been able to include.  

\section{Definition of $\kaytwo$ \citer\bass(), \citer\mbook(),
\citer\tatesym().}
Let $F$ be a (commutative) field and let $A$ be a commutative
group. An $A$-valued {\it symbol\/} on $F$ is a $\Z$-bilinear map
$s:\Fstar\times\Fstar\fleche A$ satisfying
$$
s(x,y)=0\quad\hbox{\it for all\/ $(x,y)\in\Fstar\times\Fstar$ such that\/
$x+y=1$.}
$$
A symbol $s$ is said to be trivial if $s(x,y)=0$ for all
$x,y\in\Fstar$.

Clearly, if $s$ is an $A$-valued symbol on $F$, if $f:\Fprime\fleche
F$ is a homomorphism of fields and if $a: A\fleche\Aprime$ is a
homomorphism of commutative groups, then $a\circ s\circ f$ is an
$\Aprime$-valued symbol on $\Fprime$.

\noindent
{\it Example\/}~1.  Take $F=\R$, $A=\{1, -1\}$ and set
$s_\infty(x,y)=-1$ if and only if $x<0$ and $y<0$.  Then $s_\infty$ is
a symbol.  Note that $s_\infty(x,y)=1$ if and only if the conic
$$
xR^2+yS^2=1
\numeroter \newcount\coniceq \global\coniceq=\formuleno
$$
has a solution $(R,S)$ over $\R$.  Note also that the symbol
$s_\infty$ is continuous.



\noindent
{\it Example\/} 2.  
Keep $A=\{1, -1\}$ and define an $A$-valued {\it continuous\/} symbol
on $\Q_2$ bilinearly by
$$
s_2(x,y)=(-1)^{{x-1\over2}{y-1\over2}}\ ,\
s_2(x,2)=(-1)^{{x^2-1\over2^3}}\quad(x,y\in\Z_2^\times).
$$
One has $s_2(x,y)=1$ if and only if the conic \formule\coniceq\ has a
solution over $\Q_2$.

\noindent
{\it Example\/} 3.  Let $p$ be an odd prime number and take $F=\Qp$,
$A=\Fpstar$.  Denoting by $v:\Qpstar\fleche\Z$ the valuation, one sees
that the element
$$
(-1)^{v(x)v(y)}\,{x^{v(y)} y^{-v(x)}}
\numeroter \newcount\tamesym \global\tamesym=\formuleno
$$
belongs to $\Zpstar$ for every $x,y\in\Qpstar$~; let $s_p(x,y)$ be its
image in $\Fpstar$. Then $s_p$ is a {\it continuous\/} symbol.  Further,
$h_p(x,y)=\displaystyle s_p(x,y)^{{ p-1\over 2}}$ is also a symbol,
with values in $\{1, -1\}$.  As an exercise, show that $h_p(x,y)=1$ if
and only if the conic \formule\coniceq\ has a solution $(R,S)$ over
$\qp$.


As an aside, let us mention the local-to-global principle for
$\Q$-conics (i.e,~curves of genus~0)~: for $x,y\in\Q^\times$, the
conic \formule\coniceq\ has a solution $(R,S)$ over $\Q$ if (and of
course only if) $h_p(x,y)=1$ for every odd prime $p$, $s_2(x,y)=1$,
and $s_\infty(x,y)=1$.  It so happens that if all but one of these
local conditions are satisfied, then so is the recalcitrant local
condition (cf.~th.~10).

(Let $F$ be any field and let $x,y\in\Fstar$.  The conic
\formule\coniceq\ has a solution over $F$ if and only if there exists
a homomorphism $F[R,S]/(xR^2+yS^2)\rightarrow F$ of $F$-algebras.
{\it Examples\/}~1--3 thus provide a criterion for the existence of
such a homomorphism when $F=\R,\Q_2, \qp$ ($p$ odd prime),
respectively.)

\noindent
{\it Example\/} 4.  Let $k$ be any field, take $F=k(\!(T)\!)$,
$A=\kstar$ and denote by $v:\Fstar\fleche\Z$ the valuation. For
every pair $x,y\in\Fstar$, the element \formule\tamesym\ belongs to
$k[[T]]^\times$~; let $s(x,y)$ be its image in $\kstar$. Then $s$ is a
symbol.

\noindent
{\it Example\/} 5.  More generally, let $F$ be any field endowed with
a discrete valuation $v:\Fstar\fleche\Z$, and denote the residue field
by~$k$.  We get a $\kstar$-valued symbol $s$ on $F$ by defining
$s(x,y)$ to be the image in $\kstar$ of the element \formule\tamesym~;
$s$ is called the {\it tame\/} symbol associated to $v$. {\it
Examples\/}~3~and~4 are particular cases~; {\it Examples\/}~1~and~2
are not.



Does there exist a universal symbol
$u:\Fstar\times\Fstar\fleche U_F$ on $F$\thinspace?  In other
words, does there exist a commutative group $U_F$ and a symbol $u$ on
$F$ with values in $U_F$ such that, given any $A$-valued symbol $s$ on
$F$ ($A$ being a commutative group), there exists a unique
homomorphism $f: U_F\fleche A$ of groups such that
$s=f\circ u$\thinspace?  Clearly, if such a universal symbol exists,
it is unique, up to unique isomorphism.

Nor is the {\it existence\/} of a universal symbol hard to see~: just
take $U_F$ to be the quotient of $\Fstar\otimes_{\Z}\Fstar$ by the
subgroup generated by those elements $x\otimes y$ ($x,y\in\Fstar$) for
which $x+y=1$.

\th DEFINITION 1
\enonce
Let\/ $F$ be a field.  The universal symbol on\/ $F$ is denoted\/
$\symbol(\ ,\ ):\Fstar\times\Fstar\fleche\kaytwo(F)$, or more
precisely by\/ $\symbol(\ ,\ )_F$.
\endth
The risk of confusing the {\it element\/} $\symbol(x,y)\in\kaytwo(F)$
with the {\it subset\/} $\{x,y\}\subset\Fstar$ is minor.

\smallbreak\noindent
{\it Example\/} 6.  Let $k$ be a field, take $F=k(T)$ and
$A=\kaytwo(k)$.  Write $x\in\Fstar$ as $\displaystyle
x={a_mT^m+\cdots+a_0\over b_nT^n+\cdots+b_0}$ ($a_m,b_n\in\kstar$) and
define $\displaystyle c(x)={a_m\over b_n},$ which does not depend on
the way $x$ is written as a quotient of two polynomials.  Then
$(x,y)\mapsto\symbol(c(x),{c(y)})_k$ is a $\kaytwo(k)$-valued symbol
on $F$.  Hence there is a unique homomorphism
$\kaytwo(F)\fleche\kaytwo(k)$ such that
$\symbol(x,y)_F=\symbol(c(x),{c(y)})_k$ for all $x,y\in\Fstar$. One
can show that it is a retraction of the canonical homomorphism
$\kaytwo(k)\fleche\kaytwo(F)$.

In fact, Milnor defines a graded ring $\kay(F)=\sommedir_n\kayenn(F)$,
functorially in the field $F$, as being the quotient of the
$\N$-graded tensor algebra of the $\Z$-module $\Fstar$
$$
\hbox{\bigss T}(\Fstar)=
\Z\oplus\Fstar\oplus(\Fstar\otimes\Fstar)\oplus\ldots
$$
by the homogenous ideal generated by the elements $x\otimes y$
($x,y\in\Fstar$) such that $x+y=1$.  Thus, $\kay_0(F)=\Z$ and
$\kay_1(F)$ is generated by $\{x\}$ ($x\in\Fstar$), subject to the
only conditions that $\{xy\}=\{x\}+\{y\}$ ; the map $x\mapsto\{x\}$ is
an isomorphism $\Fstar\fleche\kay_1(F)$.  We shall say very little
about the groups $\kayenn(F)$ for $n>2$.

If we are given a graded ring $A=\sommedir_nA_n$ and a homomorphism of
groups $\varphi:\Fstar\fleche A_1$ such that
$(x,y)\mapsto\varphi(x)\varphi(y)$ is an $A_2$-valued symbol on $F$,
there exists a unique homomorphism of graded rings
$\kay(F)\fleche A$ extending~$\varphi$.  Sometimes, this
homomorphism or its components are also called {\it symbols}, by
abuse of language. People go to the extent of calling elements of
$\kay(F)$ or their images in $A$ symbols.

Quillen has defined higher $\kay$-groups for arbitrary schemes ; they
do not agree for fields with those of Milnor for $n>2$.  In a certain
precise sense, Milnor's theory is the
\bquote simplest part\equote\ of Quillen's theory.

\goodbreak
\th PROPOSITION 1
\enonce
Let $F$ be a field.  In the group $\kaytwo(F)$, one has :
$$\vbox{\halign{\hfil\it #\/\rm)\ &\ $#$\hfil&\qquad$#$\hfil\cr
  i&\symbol(x,-x)=0&(x\in\Fstar),\cr
 ii&\symbol(x,y)+\symbol(y,x)=0&(x,y\in\Fstar),\cr
iii&\symbol(x,x)=\symbol(-1,x)&(x\in\Fstar).\cr
}}$$
\endth
\proof
{\it i\/}).  This is clear by bilinearity if $x=1$.  If $x\neq1$,
write $-x=(1-x)(1-\xinv)\inv$ and compute
$$
\symbol(x,-x)
=\symbol(x,1-x)-\symbol(x,1-{x^{-1}})
=\symbol({x^{-1}},1-{x^{-1}})
=0.
$$
Using {\it i\/}), we have
$$
0=\symbol(xy,-xy)
=\symbol(x,-x)+\symbol(x,y)+\symbol(y,x)+\symbol(y,-y)
=\symbol(x,y)+\symbol(y,x)
$$
so {\it ii\/}) holds.  For {\it iii\/}), we have
$\symbol(x,x)=\symbol(-1,x)+\symbol(-x,x)=\symbol(-1,x)$, again using
{\it i\/}).  Note : {\it ii\/}) implies that the projection
$\hbox{\bigss T}(\Fstar)\fleche\kay(F)$ factors through the exterior
algebra $\bigwedge(\Fstar)$ if we tensor all three with
$\Z\!\!\left[{1\over2}\right]$.
\cqfd

\section{Computations of $\kaytwo$ \citer\matsumoto(),
\citer\mbook(), \citer\basstate().}  
The fact that the group $\kaytwo(F)$ vanishes for every finite field
$F$ will have many consequences when $\kaytwo$ will be related to
other objects, like central extensions, central simple algebras,
quadratic forms, etc.

\th THEOREM 1 (Steinberg)
\enonce
For every finite field\/ $F$, the group\/ $\kaytwo(F)$ is trivial.
Consequently, all the higher\/ $\kay$-groups of a finite field are
trivial.
\endth
\proof  The multiplicative group $ F^\times $ of a finite field is always
cyclic\thinspace; let $\zeta$ be a generator.  It is sufficient to
show that $\symbol(\zeta,\zeta)=0$.

One has $2\symbol(\zeta,\zeta)=0$ by Prop.~1{\it ii\/}).  This implies
that $2\symbol(x,y)=0$ for any $x,y\in F^\times $ : writing
$x=\zeta^m$, $y=\zeta^n$, we have
$2\symbol(x,y)=2\symbol(\zeta^m,\zeta^n)=2mn\symbol(\zeta,\zeta)=0$.
So $K_2(F)$ is killed by~$2$.

If\/ $F$ is of characteristic~$2$, then $\zeta=-\zeta$ and
$\symbol(\zeta,\zeta)=\symbol(\zeta,-\zeta)=0$ by Prop.~1{\it i\/}),
so the proof is complete in this case.

Suppose that $2$ is invertible in $F$.  A counting argument shows that
the equation $\zeta x^2+\zeta y^2 =1$ has a solution $(x,y)\in
F^\times \times F^\times $ : the number of elements $\zeta x^2$ ($x\in
F$), plus the number of elements $1-\zeta y^2$ ($y\in F$), is 
$>\Card F$, so there is a pair $x,y\in F$ such that $\zeta x^2+\zeta
y^2 =1$.  Necessarily $x,y\in F^\times$, as $\zeta$ is not a square.


We then have $0=\symbol(\zeta x^2,\zeta y^2)=\symbol(\zeta,\zeta)$
since, as we have seen, $K_2(F)$ is killed by~$2$.  This completes the
proof.
\cqfd

The symbols of {\it Examples\/}~1, 2, 3 are the only {\it
continuous\/} symbols on $\R$, $\Q_2$, $\qp$ ($p$ odd) respectively
(Moore).  For more information about the $\kaytwo$ of local fields,
see the section on the uniqueness of reciprocity laws. 

%


%

Let us state the result of Tate's computation of $\kaytwo(\Q)$.  He
was inspired by the first proof of the law of quadratic reciprocity
given by the young Gauss.

{\it Example\/}~2 provides a homomorphism
$s_2:\kaytwo(\Q)\fleche\{1,-1\}$ and, similarly, {\it Example\/}~3
provides $s_p:\kaytwo(\Q)\fleche\Fpstar$ for every odd prime~$p$.
Given $x,y\in\Q^\times$, one has $s_p(x,y)=1$ for almost all~$p$.  So
we get a map into the direct sum
$$
\kaytwo(\Q)\longrightarrow
\{1,-1\}\oplus
(\Z/3\Z)^\times\oplus
(\Z/5\Z)^\times\oplus
(\Z/7\Z)^\times\oplus\ldots
$$
\vskip-20pt
\th THEOREM 2 (Tate)
\enonce
This map is an isomorphism of groups.
\endth



%


It is easily seen that the law of quadratic reciprocity is a
consequence of Tate's calculation : we have another symbol on $\Q$,
namely the symbol $s_\infty$ of {\it Example\/}~1 ; expressing it in
terms of the symbols $s_2, s_3, s_5,\ldots$ (it turns out to be the
product $s_2h_3h_5\ldots$, cf.~{\it Examples\/}~2,~3) leads to the
result.


When $k$ is a field and $F=k(T)$, the group $\kaytwo(F)$ is
canonically isomorphic to the direct sum of $\kaytwo(k)$ and the
various groups $(k[T]/\pgoth)^\times$, where $\pgoth$ runs through the
maximal ideals of $k[T]$ (cf.~{\it Example\/}~4).  Thus we get an
exact sequence
$$
\zero\fleche
\kaytwo(k)\fleche
\kaytwo(F)\fleche
\sommedir_{\pgoth} (k[T]/\pgoth)^\times\fleche
\one
\numeroter \newcount\ratfun \global\ratfun=\formuleno
$$
which admits a retraction (cf.~{\it Example\/}~6).  Applying this
result to \bquote the place at infinity\equote\ of $F|k$, one is lead to
Weil's reciprocity law.


%
\def\Fii{{\bf F}_2}
Let $A=\Fii[T]$ be the polynomial ring in one indeterminate over the
field $\Fii$, and let $A=A_0\oplus A_1\oplus A_2\oplus\ldots$ be its
gradation.  Consider the homomorphism $\phi:\R^\times\fleche A_1$
which sends $x\in\R^\times$ to $T$ if and only if $x<0$.  Then
$(x,y)\mapsto\phi(x)\phi(y)$ ($x,y\in\R^\times$) is the symbol
$s_\infty$ of {\it Example\/}~1~; it can be shown that the associated
homomorphism $\kay(\R)/2\kay(\R)\fleche\Fii[T]$ is an
isomorphism of graded rings.

Let $F$ be a global field.  It can be shown (Bass-Tate) that the real
places of~$F$ --- suppose there are $r_1$ of them --- induce an
isomorphism $\kayenn(F)\fleche(\Z/2\Z)^{r_1}$ ($n>2$).  

\section{$\kaytwo$ and universal central extensions
\citer\steinberg(), \citer\mbook().}


Recall that the {\it derived subgroup\/} $\derivedsg(G)$ of a group
$G$ is the subgroup of $G$ generated by all the commutators
$\commut(x,y)=xyx\inv y\inv$ ($x,y\in G$)~; it is normal in $G$ and
$G/\derivedsg(G)$ is the largest commutative quotient of $G$.  A
group $G$ is said to be {\it perfect\/} if $\derivedsg(G)=G$.

Recall also that an {\it extension} of a group $G$ is a pair $(G',p)$
consisting of a group $G'$ and a surjective homomorphism
$p:G'\fleche G$~; it is said to be {\it central\/} if the subgroup
$\Ker(p)\subset G'$ is central.  A {\it section\/} of an extension is
a homomorphism $s:G\fleche G'$ such that $p\circ s=\Id_G$ ; an
extension is said to {\it split\/} if it admits a section.  A
{\it morphism\/} from an extension $(G',p)$ to an extension $(G'',p')$
of $G$ is a homomorphism of groups $f:G'\fleche G''$ such that
$p=p'\circ f$.  If an extension $(G',p)$ of $G$ splits, it is
isomorphic to $(\Ker(p)\times G,\pr_2)$.

A central extension $U$ of a group $G$ is said to be
{\it universal\/} if, given any central extension $G'$ of $G$,
there exists a unique morphism $f:U\fleche G'$.  The universal
central extension, if it exists, is unique, up to unique isomorphism.

With these definitions, let us recall a few facts.  A group $G$ admits
a universal central extension if and only if $G$ is perfect.  A
central extension $U$ of a group $G$ is universal if and only if $U$
is perfect and every central extension of $U$ splits.

The universal central extension $\tilde G$ of a perfect group $G$ can
be identified with the homology group $\ache_2(G,\Z)$ (whose
definition we do not recall) ; $\Ker(\tilde G\fleche G)$ is the
central subgroup of $\tilde G$.  Under the identification
$\ache^2(G,\ache_2(G,\Z))\fleche\End_{\Z}(\ache_2(G,\Z))$, the
extension $\tilde G$ corresponds to
$\Id:\ache_2(G,\Z)\fleche\ache_2(G,\Z)$.  Also, a central extension
$(G',p')$ of $G$ is universal if and only if $\ache_1(G',\Z)=\zero$
and $\ache_2(G',\Z)=\zero$.

Some people think of $G/\derivedsg(G)=\ache_1(G, \Z)$ as the $\pi_0(G)$ of a
discrete group $G$ and call it \bquote connected\equote\ if $\pi_0(G)=\one$ ;
$\tilde G$ is then its \bquote universal cover\equote\ and the kernel of
$\tilde G\fleche G$ is the \bquote fundamental group\equote\ $\pi_1(G)$
of~$G$.


Now let $F$ be a field and $n>2$ an integer.  The group $\SL_n(F)$ is
perfect --- with three exceptions : when $(n,\Card F)$ equals $(3,2)$,
$(3,4)$ or $(4,2)$ ; we do not consider them in what follows.  Let
$(\tilde G,p)$ be its universal central extension.  For
$x,y\in\Fstar$, choose $x_{12}, y_{13}\in\tilde G$ such that
$$
p(x_{12})=\diag(x,\xinv,1,1,\ldots)\ \ \hbox{and}\ \
p(y_{13})=\diag(y,1,y\inv,1,\ldots)
$$ in $\SL_n(F)$.  Then the commutator $\commut(x_{12},y_{13})$ ---
which depends only on $x$~and $y$ --- belongs to
$\Ker(p)=\ache_2(\SL_n(F),\Z)$ and defines a bilinear map
$$
(x,y)\mapsto\commut(x_{12},y_{13})
:\Fstar\times\Fstar\fleche\ache_2(\SL_n(F),\Z).
\numeroter \newcount\ucesym \global\ucesym=\formuleno
$$
\vskip-20pt
\th THEOREM 3 (Matsumoto, Steinberg)
\enonce
Apart from the three exceptional\/ $(n,F)$ mentioned above, the map
\formule\ucesym\ is a symbol inducing an isomorphism\/
$
\kaytwo(F)\fleche\ache_2(\SL_n(F),\Z).
$
\endth


Note that, in view of the triviality of the $\kaytwo$ of a finite
field, this theorem implies that the {\it discrete\/} groups
$\SL_n(F)$ ($n>2$, $F$ finite) are \bquote simply connected\equote ---
leaving aside the three exceptions noted above, which are not \bquote
connected\equote.

In fact, $\kaytwo(F)$ is the \bquote fundamental group\equote\ of
$G(F)$ for every simple, simply connected split algebraic $F$-group
$G$, not just $G=\SL_n(F)$.  The only exceptions are $\Sp_n(F)$, for
which $\kaytwo(F)$ is a quotient of the
\bquote fundamental group\equote, with kernel $\Z$.


\goodbreak



\section{$\kaytwo$ and quadratic forms \citer\milnor(), \citer\kahn(),
\citer\morel().}

Let $F$ be a field in which $2$ is invertible.  Recall that a
(regular) quadratic space over $F$ is a (finite-dimensional) vector
$F$-space $V$ endowed with a
regular symmetric bilinear form $b$.  It is possible to choose a
basis for $V$ and $a_1,\ldots,a_n\in\Fstar$ such that one has
$$
b(\xi,\xi)=a_1\xi_1^2+a_2\xi_2^2+\cdots+a_n\xi_n^2
\numeroter \newcount\qform \global\qform=\formuleno
$$
for all $\xi=(\xi_1,\xi_2,...,\xi_n)$ ($n=\dim V$) in $V$.  The isometry
classes of quadratic $F$-spaces form a monoid with orthogonal direct
sum as the law of addition ;
this monoid is integral (Witt). The corresponding group of differences
$\witthat(F)$ is called the Grothendieck group of $F$.  The class
of the form \formule\qform\ in $\witthat(F)$ is denoted
$\qfclass(a_1,a_2,\ldots,a_n)
=\qfclass(a_1)+\qfclass(a_2)+\ldots+\qfclass(a_n)$.  Tensor product
over $F$ makes $\witthat(F)$ into a ring whose multiplication is
characterised by $\qfclass(x)\qfclass(y)=\qfclass(xy)$ and for which
$\qfclass(1)$ is the neutral element.  Taking dimensions gives a
canonical surjection of rings $\dim:\witthat(F)\fleche\Z$. Denote by
$\hat I$ the kernel (the augmentation ideal).

Let $h=\qfclass(1,-1)=\qfclass(1)+\qfclass(-1)$ be the class of the
hyperbolic plane ; it corresponds to the quadratic space $L\oplus
\Hom_F(L,F)$, $(\xi,\xi^*)\mapsto\xi^*(\xi)$, where $L$ is a vector
$F$-line.  The subgroup $H$ generated by $h$ is an ideal in
$\witthat(F)$ (a quadratic space $\qfclass(a_1,a_2,\ldots,a_n)$
\bquote represents~$0$\equote\ if and only if
$\qfclass(a_1,a_2,\ldots,a_n)=\qfclass(b_1,b_2,\ldots,b_{n-2})+h$ for
some quadratic space $\qfclass(b_1,b_2,\ldots,b_{n-1})$) ; one has
$H\cap\hat I=\zero$.  The quotient $\witt(F)=\witthat(F)/H$ is called
the Witt ring of $F$ ; the image $I=\hat I/(H\cap\hat I)=\hat I$ of
the ideal $\hat I$ is maximal in $\witt(F)$ with $\Z/2\Z$ as the
quotient.  We will be mainly interested in the graded
$\F_{\!2}$-algebra associated to the filtered ring
$$
\ldots\subset I^3\subset I^2\subset I\subset\witt(F).
$$
There is a homomorphism $s_1:\Fstar\fleche I/I^2$ given by
$s_1(x)=\qfclass(x)-\qfclass(1)$.
\th LEMMA 1 (Milnor, 1970)
\enonce
The map\/ $s_2(x,y)=s_1(x)s_1(y)$
is an\/ $(I^2\!/I^3)$-valued symbol on\/ $F$ ; the corresponding
homomorphism from\/ $\kaytwo(F)$ is trivial on\/ $2\kaytwo(F)$.
\endth
There is therefore a unique homomorphism of graded $\F_{\!2}$-algebras
$$
s:\kay(F)/2\kay(F)\longrightarrow
\Z/2\Z\oplus I/I^2\oplus I^2\!/I^3\oplus\cdots
\numeroter \newcount\wittiso \global\wittiso=\formuleno
$$
extending $s_1$ and $s_2$.  It is easy to see that the discriminant
provides an inverse $w_1:I/I^2\fleche\Fstar\!/\Fstar{}^2$ of $s_1$
and the Hasse invariant an inverse
$w_2:I^2\!/I^3\fleche\kaytwo(F)/2\kaytwo(F)$ of $s_2$.

Milnor showed that the map \formule\wittiso\ is always surjective and
conjectured (1970) that it is bijective for all fields $F$ (in which
$2$ is invertible).  He proved the bijectivity for global fields.

The conjecture was finally proved by Orlov, Vishik \& Voevodsky, in a
preprint available on the {\tt arXiv}.

\th THEOREM 4 (Orlov, Vishik, Voevodsky, 1996)
\enonce
The map \formule\wittiso\ is an isomorphism of graded
$\F_{\!2}$-algebras.
\endth

\centerline{***}

Continuing to assume that $2$ is invertible in the field $F$, let
$\Fsep$ be a {\it separable algebraic closure\/} of $F$ and
$\Gamma=\Gal(\Fsep|F)$ the (profinite) group of $F$-automorphism of
$\Fsep$.  Consider the exact sequence
$$
\one\fleche
\{1,-1\}\fleche
\Fsepstar\droite{()^2}
\Fsepstar\fleche\one
$$
of discrete $\Gamma$-modules.  The associated long cohomology sequence
furnishes --- upon identifying the $\Gamma$-module
$\{1,-1\}={}_2\Fsepstar$ with $\Z/2\Z$  --- an injection
$\delta_1:\Fstar\!/\Fstar{}^2\fleche\ache^1(\Gamma,\Z/2\Z)$ which
is an isomorphism by Hilbert's theorem~90 :
$\ache^1(\Gamma,\Fsepstar)=\zero$.

\th LEMMA 2 (Tate, 1970)
\enonce
The map\/ $\delta_2(x,y)=\delta_1(x)\cupprod\delta_1(y)$ is an\/
$\ache^2(\Gamma,\Z/2\Z)$-valued symbol on\/ $F$ ; the corresponding
homomorphism from\/ $\kaytwo(F)$ is trivial on\/ $2\kaytwo(F)$.
\endth
There is therefore a unique homomorphism of graded $\F_{\!2}$-algebras
$$
\delta:\kay(F)/2\kay(F)\longrightarrow
\Z/2\Z\oplus\ache^1(\Gamma,\Z/2\Z)\oplus\ache^2(\Gamma,\Z/2\Z)\oplus\cdots
\numeroter \newcount\galiso \global\galiso=\formuleno
$$
extending $\delta_1$ and $\delta_2$.  We have seen that $\delta_1$ is
an isomorphism ; Merkurjev (1981) proved that $\delta_2$ is an
isomorphism.  Milnor (1970) proved that $\delta$ is an isomorphism
when $F$ is a finite, a local, a global or a maximally ordered field.
His conjecture --- that this is an isomorphism for all fields $F$ (in
which $2$ is invertible) --- was finally proved by Voevodsky in 1996
and got him a medal.  Nearly a decade earlier, Rost and Merkurjev \&
Suslin had proved that $\delta_3$ is an isomorphism.
\th THEOREM 5 (Voevodsky, 1996)
\enonce
The map\/ \formule\galiso\ is an isomorphism of graded
$\F_{\!2}$-algebras.
\endth
As a consequence of theorems~4 and~5, the graded ring associated to
the filtered ring $\witt(F)$ is canonically isomorphic to the
cohomology ring of the $\Gamma$-module $\Z/2\Z$.  To see how far this
result goes, note that previous attempts to construct maps between the
two --- let alone showing that they are isomorphisms --- had been
successful only in low degrees.

\section{$\kaytwo$ and central simple algebras \citer\mbook(),
\citer\soule(), \citer\tatekii().}

Let $F$ be a (commutative) field.  Recall that a (finite-dimensional,
associative, unital) $F$-algebra $A$ is called {\it simple\/} if the
only bilateral ideals of $A$ are $\zero$ and $A$.  An $F$-algebra $A$
is called {\it central\/} if $F$ is precisely the centre of $A$.
Every central simple $F$-algebra $A$ is isomorphic to the matrix
algebra $\M_n(D)$ of a (skew) field $D$ over $F$ (Wedderburn) ; the
pair $(n,D)$ is uniquely determined by $A$, up to isomorphism.  Two
such algebras $A$, $A'$ are called {\it similar\/} if the
corresponding (skew) fields are isomorphic.  Similarity classes of
central simple $F$-algebras form a group $\br(F)$ ({\it the Brauer
group of\/ $F$}) with tensor product as the law of multiplication.  It
is a torsion group, as a restriction-corestriction argument shows.


The group $\br(F)$ can also be viewed as the group of $F$-isomorphism
classes of $F$-algebras which become isomorphic, over a separable
algebraic closure $\Fbar$ of $F$, to $\M_n(\Fbar)$ (some $n$), with
tensor product of algebras providing the group law.  It can also be
viewed as the group of $F$-varieties which are $\Fbar$-isomorphic to
$\P_n$ (some $n$).


Now let $n>0$ be an integer and assume that $F$ contains a {\it
primitive\/} $n^{\rm th}$ root $\zeta$ of\/ $1$ (i.e.\ $\zeta\in\Fstar$
is of order $n$ ; thus $n$ is invertible in $F$).  For $a,b\in\Fstar$,
consider the $F$-algebra $A_\zeta(a,b)$ with the presentation
$$
x^n=a\ ;\quad y^n=b\ ;\quad xy=\zeta yx.
$$
It is a central simple $F$-algebra --- called a {\it cyclic\/} algebra
in general and a {\it quaternion\/} algebra when $n=2$, since it was
so called in the case $F=\R$, $\zeta=-1$, $a=-1$, $b=-1$ by its
discoverer Hamilton --- whose class $s_\zeta(a,b)\in\br(F)$ is killed
by $n$, i.e.~lies in the $n$-torsion ${}_n\br(F)$.

(Before going on, observe that a quaternion $\Qp$-algebra defined by
$a,b\in\Qp^\times$ is the matrix algebra precisely when $h_p(a,b)=1$
({\it Example\/}~3) for an odd prime $p$, when $s_2(a,b)=1$ ({\it
Example\/}~2) for $p=2$, and when $s_\infty(a,b)=1$ ({\it
Example\/}~1) for $\Qp=\R$.  The local-to-global principle for
quaternion $\Q$-algebras $A$, essentially the same as the one for
conics, says that $A$ is a matrix algebra if (and of course only if)
$A\otimes\Qp$ is a matrix algebra for every place $p$ of $\Q$.  As
before, it is in fact sufficient to demand that they be matrix
algebras at all places, with one possible exception; the ``\thinspace
exception\thinspace'' is then not an exception.  There is a
local-to-global principle for central simple algebras over every
global field, i.e.~over finite extensions of\/ $\Q$ or of\/ $\F_p(T)$
($p$ prime).)

\th LEMMA 3 (Tate, 1970)
\enonce
The map\/ $(a,b)\mapsto s_\zeta(a,b)$ is an\/ $\,{}_n\br(F)$-valued symbol
on\/ $F$.
\endth

\th THEOREM 6 (Merkurjev \& Suslin, 1982)
\enonce
The associated map is
an isomorphism\/ $\kaytwo(F)/n\kaytwo(F)\fleche{}_n\br(F)$.
\endth

Thus every central simple algebra whose class is killed by $n$ is
similar to a product of cyclic algebras in the presence of
$n^{\rm th}$ roots of $1$.

The choice $\zeta$ of a primitive $n^{\rm th}$ root of $1$ in $F$
allows us to identify $\,{}_n\br(F)$ with $\ache^2(\Gamma,\Z/n\Z)$.
The lemma implies that there is a unique homomorphism of graded
$(\Z/n\Z)$-algebras
$$
\delta_\zeta:\kay(F)/n\kay(F)\longrightarrow
\Z/n\Z\oplus\ache^1(\Gamma,\Z/n\Z)\oplus\ache^2(\Gamma,\Z/n\Z)\oplus\cdots
\numeroter \newcount\delze \global\delze=\formuleno
$$
which restricts to the symbol $s_\zeta$ on $\Fstar\times\Fstar$.  But
let turn to what happens when a primitive $n^{\rm th}$ root of $1$ may
not be available in $F$.
\smallskip
\centerline{***}
\goodbreak

More generally, without assuming the existence of a primitive
$n^{\rm th}$ root of $1$ in $F$ but merely that $n$ is invertible in
$F$, we have an exact sequence
$$
\one\fleche
\Z/n\Z(1)\fleche
\Fsepstar\droite{()^n}
\Fsepstar\fleche\one
$$
of discrete $\Gamma$-modules, where $\Z/n\Z(1)={}_n\Fsepstar$ is the
group of $n^{\rm th}$ roots of $1$ in $\Fsepstar$ with its natural
$\Gamma$-action.  The associated long exact cohomology sequence and
Hilbert's theorem~90 furnish --- as we have seen before in the case
$n=2$ --- an isomorphism
$\delta_1:\Fstar\!/\Fstar{}^n\fleche\ache^1(\Gamma,\Z/n\Z(1))$.
(The next few terms of the same cohomology sequence furnish a
canonical isomorphism
$\ache^2(\Gamma,\Z/n\Z(1))\fleche{}_n\br(F)$, a group which also
classifies $F$-varieties potentially isomorphic to $\P_{n-1}$).  Cup
product on cohomology
$$
\cupprod\;:\ache^r(\Gamma,\Z/n\Z(r))
\times\ache^s(\Gamma,\Z/n\Z(s))\longrightarrow
\ache^{r+s}(\Gamma,\Z/n\Z(r+s))
$$
then provides a bilinear map
$
\delta_2:\Fstar\!/\Fstar{}^n\times\Fstar\!/\Fstar{}^n\fleche
\ache^2(\Gamma,\Z/n\Z(2)).
$
\th LEMMA 4 (Tate, 1970)
\enonce
The map\/ $\delta_2(x,y)=\delta_1(x)\cupprod\delta_1(y)$ is a
symbol on\/ $F$.
\endth
Choosing a primitive $n^{\rm th}$ root $\zeta$ of $1$ when there is
one in $F$, and using it to identify the groups ${}_n\br(F)$,
$\ache^2(\Gamma,\Z/n\Z)$ and $\ache^2(\Gamma,\Z/n\Z(2))$, the symbol
$\delta_2$ turns out to be the same as the symbol $s_\zeta$ of
lemma~3.  The map \formule\delze\ associated to $s_\zeta$ is thus the
same in this case as the one in the :

\th CONJECTURE (Bloch \& Kato, 1986)
\enonce
The associated homomorphism of graded $(\Z/n\Z)$-algebras
$$
\delta:\kay(F)/n\kay(F)\longrightarrow
\Z/n\Z\oplus\ache^1(\Gamma,\Z/n\Z(1))\oplus
\ache^2(\Gamma,\Z/n\Z(2))\oplus\cdots
\numeroter 
$$
is an isomorphism for all fields\/ $F$ in which\/ $n$ is invertible.
\endth
\newcount\bkgalois \global\bkgalois=\formuleno


The main theorem (cf. th.~6) of Merkurjev \& Suslin (1982) says that
the map
$\delta_2:\kaytwo(F)/n\kaytwo(F)\fleche\ache^2(\Gamma,\Z/n\Z(2))$
is always an isomorphism ; Tate had proved this earlier (1976) for
global fields.  Bloch, Gabber \& Kato prove this conjecture when $F$
is a field of characteristic $0$ endowed with a henselian discrete
valuation of residual characteristic $p\neq0$ and $n$ is a power of $p$.

The Bloch-Kato conjecture --- which remains a major open problem in
its full generality --- makes the remarkable prediction that the
graded algebra $\sommedir_r\ache^r(\Gamma,\Z/n\Z(r))$ is generated by
elements of degree~1.  Galois groups should thus be very special among
profinite groups in this respect.
\goodbreak

\section{$\kaytwo$ and differential forms \citer\colliot(), \citer\bgk().}
Let $A$ be a commutative ring and $B$ a commutative $A$-algebra.
Recall that an $A$-derivation on $B$ is a pair $(\phi,M)$ consisting
of a $B$-module $M$ and an $A$-linear map $\phi:B\fleche M$ such
that
$$
\phi(xy)=x\phi(y)+\phi(x)y\quad\hbox{\it for all\/}\ x,y\in B.
$$
A morphism of $A$-derivations $(\phi,M)$, $(\phi',M')$ on $B$ is an
$B$-linear map $f:M\fleche M'$ such that $\phi'=f\circ\phi$.  Does
there exist a universal $A$-derivation $(d_B,\Omega^1_{B|A})$ on $B$ ?
Clearly it would then be unique, up to unique isomorphism.  As for the
existence, we just need to take the $B$-module generated by $d_B(x)$
($x\in B$) subject to the relations
$$\vbox{\halign{\hfil\it#\/\rm)\ \ &#\hfil\cr
i&$d_B(\alpha x)=\alpha\,d_B(x)$ for all
$\alpha\in A,x\in B$,\cr
ii&$d_B(x+y)=d_B(x)+d_B(y)$ for all
$x\in B,y\in B$,\cr
iii&$d_B(xy)=d_B(x)\,y\,+\,x\,d_B(y)$ for all
$x\in B,y\in B$.\cr
}}
$$
Write $\Omega^0_{B|A}=B$ and put
$\Omega^n_{B|A}=\bigwedge^n_{B}\Omega^1_{B|A}$ for $n>0$.  Recall that
there is a unique system of $A$-linear maps
$d_B:\Omega^n_{B|A}\fleche\Omega^{n+1}_{B|A}$ extending
$d_B:B\fleche\Omega^{1}_{B|A}$, verifying $d_B\circ d_B=0$, and
such that
$$
d_B(\omega\wedge\omega')=
d_B(\omega)\wedge\omega'+(-1)^{r+r'}\omega\wedge d_B(\omega')\quad
(\omega\in\Omega^r_{B|A},\omega'\in\Omega^{r'}_{B|A}).
$$

Let $F$ be a field.  In what follows, we shall take $A=\Z$, $B=F$
and simplify notations to $d$, $\Omega^n_F$.  The groups $\Omega^n_F$
are generated by the elements
$$
x\derlog(y_1)\wedge\cdots\wedge\derlog(y_n)
\quad(x\in F,\ y_1,\ldots,y_n\in\Fstar).
\numeroter \newcount\ungen \global\ungen=\formuleno
$$
i.e.~the map
$F\otimes\Fstar\otimes\cdots\otimes\Fstar\fleche\Omega^n_F$ which
sends $x\otimes y_1 \otimes\cdots\otimes y_n$ to \formule\ungen\ is
surjective.  We have a homomorphism
$\dlog_1:\Fstar\fleche\Omega^1_F$,
$\;\displaystyle\dlog_1(x)=\derlog(x)$, of groups.
\th LEMMA 5 (Tate, 1970)
\enonce
The map\/
$\displaystyle
\dlog_2(x,y)=\derlog(x)\wedge\derlog(y)
$
is an\/ $\Omega^2_F$-valued symbol on\/ $F$.
\endth
There is therefore a unique map
$\dlog:\kay(F)\fleche\sommedir_n\Omega^n_F$ of graded rings which
restricts to $\dlog_2$ on $\kaytwo(F)$ ; it restricts to $\dlog_1$ on
$\Fstar$.

Now suppose that $F$ is of characteristic $p\neq0$ and denote by $F^p$
the subfield of $F$ which is the image of the homomorphism
$()^p:F\fleche F$.  Denote by $Z^n_F\subset\Omega^n_F$ the kernel
of the differential $d:\Omega^n_F\fleche\Omega^{n+1}_F$ and by
$B^n_F\subset Z^n_F$ the image of the differential
$d:\Omega^{n-1}_F\fleche\Omega^n_F$.

As $d(\alpha^p\omega)=\alpha^p\,d(\omega)$ for all $\alpha\in F$,
$\omega\in\Omega^n_F$, the differential $d$ is $F^p$-linear.  Thus,
$B^n_F$ and $Z^n_F$ are vector sub-$F^p$-spaces.  We shall consider
them as vector $F$-spaces via the isomorphism $()^p:F\fleche F^p$.  

Recall that there is a (unique) system of isomorphisms
$\gamma:\Omega^n_F\fleche Z^n_F/B^n_F$ of vector $F$-spaces such
that
$$
\gamma\left(x\derlog(y_1)\wedge\cdots\wedge\derlog(y_n)\right)=
x^p\derlog(y_1)\wedge\cdots\wedge\derlog(y_n)+B^n_F
$$
for the system of generators \formule\ungen\ of $\Omega^n_F$
(Cartier).  We have an $\F_{\!p}$-linear map
$\wp=\gamma-\Id:\Omega^n_F\fleche\Omega^n_F/B^n_F$ defined by
$\wp(\omega)=\gamma(\omega)-\omega+B^n_F$ ; in degree~$0$, it is the
endomorphism $x\mapsto x^p-x$ of the additive group $F$.  Explicitly,
for the system of generators \formule\ungen\ of the group
$\Omega^n_F$, one has
$$
\wp\left(x\derlog(y_1)\wedge\cdots\wedge\derlog(y_n)\right)=
(x^p-x)\derlog(y_1)\wedge\cdots\wedge\derlog(y_n)+B^n_F.
$$
Let $\nu(n)_F$ be kernel of this map.  It is easy to see that
$\wp\circ\dlog_n=0$, i.e.~$\dlog_n$ maps $\kayenn(F)$ into
$\nu(n)_F=\Ker(\wp)$.  As $\dlog_n(p\kayenn(F))=\zero$, we
have a map of graded $\F_{\!p}$-algebras
$$
\dlog:\kay(F)/p\kay(F)\fleche
\nu(0)_F\oplus\nu(1)_F\oplus\nu(2)_F\oplus\cdots
\numeroter \newcount\bkg \global\bkg=\formuleno
$$
\vskip-20pt
\th THEOREM 7 (Bloch, Gabber \& Kato, 1986)  
\enonce
The map\/ \formule\bkg\ is an isomorphism for every field\/ $F$ of
characteristic\/ $p\neq0$.
\endth

In degree~$0$ we get the isomorphism
$\Z/p\Z\fleche\F_{\!p}=\Ker(\wp:F\fleche F)$.  There is also
an interpretation of the quotients $\kay(F)/p^r\kay(F)$ but it is
more complicated to state.


\centerline{***}

Let us briefly mention the case $p=2$ where the notions of (regular)
bilinear form and (regular) quadratic space are distinct, so we get
two objects : a filtered ring $\witt(F)$ and a $\witt(F)$-module
$\quadsp(F)$.  A theorem of Kato (1982) says that the kernel
$\nu(n)_F$ of $\wp:\Omega^n_F\fleche\Omega^n_F/B^n_F$ is
isomorphic to $I^n\!/I^{n+1}$ and the cokernel to
$I^n\quadsp(F)/I^{n+1}\quadsp(F)$.

\goodbreak
\section{$\kaytwo$ and abelian extensions \citer\raskind().}
Let $F$ be a field, $\Fab$ the maximal {\it abelian\/} extension of
$F$ and $\Gamma=\Gal(\Fab|F)$ the (profinite, commutative) group of
$F$-automorphisms of $\Fab$.

Assume that $F$ is finite.  We know that then there is a natural map
$\rho:\kay_0(F)=\Z\fleche\Gamma$, namely $n\mapsto\phi^n$ where
$\phi$ is the automorphism $x\mapsto x^q$, $q=\Card F$.  The image of
$\rho$ is dense in $\Gamma$.

Next, let $F$ be a local field, i.e. a finite extension of
$\qp$ ($p$ prime) or a field isomorphic to $k(\!(T)\!)$ where $k$ is a
finite field.  Then there is a natural
$\rho:\kay_1(F)=\Fstar\fleche\Gamma$ --- which is less easy to
describe --- whose image is dense in $\Gamma$.  This map and its
properties form the essential content of the theory of abelian
extensions of local fields.

What we have just seen are the $0$- and $1$-dimensional versions of a
general theory of $n$-dimensional local fields.  Such a field $F$
($n>1$) is complete with respect to a discrete valuation whose residue
field is an $(n-1)$-dimensional local field.

\th THEOREM 8 (Kato, S.~Saito, 1986)
\enonce
Let\/ $F$ be an $n$-dimensional local field.  There is a natural
homomorphism\/ $\rho:\kayenn(F)\fleche\Gamma=\Gal(\Fab|F)$ whose
image is dense.
\endth

In addition, this \bquote reciprocity map\equote\ is compatible with
norms from finite abelian extensions, a notion which we do not pursue
here.

Thus for a $2$-dimensional local field $F$, the group $\Gal(\Fab|F)$
is generated --- as a profinite commutative group --- by elements of
the form $\rho(\symbol(x,y))$ ($x,y\in\Fstar$), and all relations
among these elements are consequences of bilinearity and
\bquote $\symbol(x,y)=1$ whenever $x+y=1$\equote.  Examples of such
fields include $\qp(\!(T)\!)$ ($p$ prime) and $k(\!(X)\!)(\!(Y)\!)$,
where $k$ is a finite field.


%

\goodbreak
\section{$\kaytwo$ and abelian coverings of curves 
\citer\blochcft(), \citer\saito(), \citer\yoshida().}

Let $F$ be a finite extension of\/ $\Qp$ or of\/ $\Fp(\!(T)\!)$,
i.e.~a field complete with respect to a discrete valuation with finite
residue field, and let $C$ be a smooth proper absolutely connected
$F$-curve.  For $x\in C^{(1)}$ a closed point of $C$, we have the tame
symbol $\kaytwo(F(C))\to F(x)^\times$ at the place $x$ of the function
field $F(C)$ ({\it Example\/} 5).  Define
$$
\vv_1(C)=\Coker\!\big(K_2(F(C))\to\bigoplus_{x\in C^{(1)}} F(x)^\times\big).
$$
For every closed point $x\in C^{(1)}$, there is a norm map
$F(x)^\times\to\Fstar$, inducing a map $\vv_1(C)\to\Fstar$.  Let
$\vv_1(C)_0=\Ker(\vv_1(C)\to\Fstar)$, a group introduced by S.~Bloch.

Let $\varpi_1(C)$ be the profinite group classifying abelian {\'e}tale
coverings of $C$~; it is the maximal commutative quotient of the
{\'e}tale fundamental group $\pi_1(C)$ of $C$.  There is a natural
surjection $\varpi_1(C)\to\varpi_1(F)$~; denote by $\varpi_1(C)_0$ the
kernel, called the ``\thinspace the geometric part\thinspace''.  By
general results (i.e.~valid for any smooth proper absolutely connected
$F$-variety), the torsion subgroup of $\varpi_1(C)_0$ is finite, and
the quotient modulo the torsion subgroup is $\hat{\Z}^r$, where $r$ is
the toric rank of the jacobian variety (of the albanese variety, in
general) $J$ of~$C$, i.e.~the dimension of the maximal split torus in
the closed fibre of the N{\'e}ron model of $J$.

Generalising from the classical case of reciprocity maps for
$1$-dimensional schemes like (the spectrum of) the ring of integers in
a number field, one constructs ``\thinspace reciprocity
maps\thinspace'', $\sigma:\vv_1(C)\to\varpi_1(C)$,
$\sigma_0:\vv_1(C)_0\to\varpi_1(C)_0$.  We have a commutative diagram
with exact rows
$$\diagram{
\zero&\to&\vv_1(C)_0&\to&\vv_1(C)&\to&\Fstar\cr
&&\vfl{\sigma_0}{}{5mm}&&\vfl{\sigma}{}{5mm}&&\vfl{\rho}{}{5mm}\cr
\zero&\to&\varpi_1(C)_0&\to&\varpi_1(C)&\to&\varpi_1(F)&\to&\zero\cr
}
$$
in which $\rho$ is the reciprocity map for the local field $F$,
inducing an isomorphism upon completion of\/ $\Fstar$.

\th THEOREM 9 (Bloch 1981, S.~Saito 1985, Yoshida 2003)
\enonce
The kernel of the reciprocity map\/ $\sigma_0$ is the maximal
divisible subgroup of\/ $\vv_1(C)_0$, and the image is the torsion
subgroup of\/ $\varpi_1(C)_0$.  The kernel of\/ $\sigma$ is the
maximal divisible subgroup of\/ $\vv_1(C)$, and the quotient of\/
$\varpi_1(C)$ by the closure of\/ $\mathop{\rm Im}(\sigma)$ is\/
$\hat{\Z}^r$, where\/ $r$ is the toric rank of the jacobian of\/ $C$.
\endth
The idea is that we have described the group $\varpi_1(C)_0$ purely in
terms of the group $\vv_1(C)_0$ constructed from the function field
$F(C)$, just as we describe the group $\Gal(H|K)$ of $K$-automorphisms
of the maximal unramified abelian extension $H$ of a number field $K$
in terms of a group constructed from $K$, namely the ideal class
group.

%
%

\goodbreak
\section{$\kaytwo$ and uniqueness of reciprocity laws \citer\moore(),
\citer\chasewater().}


Let us briefly recall the definition of the norm residue symbol
(sometimes called the Hilbert symbol).  Let $F$ be a local field other
than $\C$ and let $\mu(F)={}_\infty\Fstar$ be the group of roots of
$1$ in $F$ ; it is cyclic, and its order $m$ is invertible in $F$.
Consider the extension $L=F(\!\root m\of\Fstar\,)$ of $F$ --- it is
the maximal abelian extension of exponent $m$.  It is of finite degree
as the closed subgroup $\Fstar{}^m\subset\Fstar$ is of finite
index. ({\it Exercise\/} : compute this index.)

The quotient $\Fstar\!/\Fstar{}^m$ admits {\it two\/} descriptions.
On the one hand, by the theory of abelian extensions of local fields,
it is canonically isomorphic to the group $G=\Gal(L|F)$ of
$F$-automorphisms of $L$.  On the other hand, we get an isomorphism
$\delta:\Fstar\!/\Fstar{}^m\fleche\ache^1(G,\mu(F))$ from the short
exact sequence
$$
\one\fleche
\mu(F)\fleche
\Lstar\droite{()^m}
\Lstar\fleche
\one
\advance\abovedisplayskip by4pt
$$
of discrete $G$-modules.  As the $G$-action on $\mu(F)$ is trivial, we
have $\ache^1(G,\mu(F))=\Hom(G,\mu(F))$.  Thus we get an isomorphism
$$\Fstar\!/\Fstar{}^m\fleche\Hom(\Fstar\!/\Fstar{}^m,\mu(F)),$$ i.e. a
perfect duality of $\Z/m\Z$-modules.  The corresponding bilinear map
$\normressym(\ ,\ )_F:\Fstar\times\Fstar\fleche\mu(F)$ happens to be a
symbol --- the norm residue symbol.  It is the universal continuous
symbol on $F$ (Moore).  Further, it can be shown that the natural
surjective map $\kaytwo(F)\fleche\mu(F)$ admits a section and that
its kernel is a uniquely divisible group (Tate, Merkurjev).  Of
course, for any divisor $d$ of $m$ we get a continuous
${}_d\mu(F)$-valued symbol on $F$ by raising the norm residue symbol
to the power $n/d$.

When $F=\R$ (resp.~$\Q_2$), we have $m=2$ and the symbol $s_\infty$
(resp.~$s_2$) of {\it Example\/}~1 (resp.~{\it Example\/}~2) is the
norm residue symbol.  When $F=\qp$ ($p$ odd), we have $m=p-1$ and the
symbol $s_p$ of {\it Example\/}~3 can be viewed as the norm residue
symbol --- up to the canonical isomorphism $\mu(F)\fleche\Fpstar$.
{\it Example\/}~4 includes the case of local function fields.


Now let $F$ be a global field, $\mu(F)\subset\Fstar$ the group of
roots of~$1$ in $F$ and, for every real or ultrametric place $v$ of
$F$, let $\mu(F_v)$ be the group of roots of~$1$ in the local field
$F_v$.  Put $m=\Card\mu(F)$, $m_v=\Card\mu(F_v)$ ; $m$ divides $m_v$
for every $v$.  (Imaginary places of $F$ play no role in what
follows~: $\C^\times$ is connected and therefore the universal
continuous symbol on $\C$ is trivial).

Consider $\sommedir_v\mu(F_v)$ (where $v$ runs over the real or
ultrametric places of $F$).  We have a natural homomorphism from this
direct sum to $\mu(F)$ which on the $v$-th component is the map
$\displaystyle\zeta\mapsto\zeta^{m_v\over m_{\phantom{v}}}$.  Also,
for $x,y\in\Fstar$, we have $\normressym(a,b)_{{}_{\phantom{v}}F_v}=1$
for almost all $v$.  We thus get a sequence
$$
\kaytwo(F)\droite{\lambda}
\sommedir_v\mu(F_v)\droite{}\mu(F)\fleche\one.
\numeroter \newcount\exprec \global\exprec=\formuleno
$$
The \bquote explicit\equote\ reciprocity law says that this sequence
is a complex.  Tate's computation (th.~2) amounts to saying that
\formule\exprec\ is exact for $F=\Q$ and that $\Ker(\lambda)=\zero$.
We have seen that, in the notation of {\it Examples\/}~1--3, Hilbert's
product formula $s_\infty s_2h_3h_5\ldots=1$, which is equivalent to
the law of quadratic reciprocity, is a particular case.  In general,
one has \bquote uniqueness of reciprocity laws\equote\ in the
following sense :

\th THEOREM 10 (Moore, 1968)
\enonce
The sequence \formule\exprec\ is exact for every global field\/ $F$.
\endth
When $F$ is a function field of characteristic $p$, Tate proved that
$\Ker(\lambda)$ is a finite group of order prime to $p$.  For number
fields, the finiteness of $\Ker(\lambda)$ was proved by Brumer in the
totally real abelian case and by Garland in general, as a consequence
of his theorem --- whose proof uses riemannian geometry and harmonic
forms --- about the vanishing of $\ache^2(\SL_n(\ogoth),\R)$ ($n>6$)
for the ring of integers $\ogoth$ of $F$.  Now these finiteness
results are corollaries of the general results of Quillen.


\section{$\kaytwo$ and special values of $\zeta$-functions
\citer\bass(), \citer\licht(), \citer\mazurwiles(), \citer\rubin().}

Let $F$ be a global function field over a (finite) field $k$.  Let us
rewrite the sequence \formule\exprec\ in terms of the multiplicative
groups $k_v^\times$ of the residue fields $k_v$ at the various places
$v$ --- necessarily ultrametric --- of $F$.  At every $v$, there is a
canonical isomorphism $\mu(F_v)\fleche k_v^\times$ ; there is also an
isomorphism $\mu(F)\fleche k^\times$.  So the exact sequence
\formule\exprec\ becomes
$$
\zero\fleche
\Ker(\lambda)\fleche
\kaytwo(F)\droite{\lambda}
\sommedir_vk_v^\times\fleche
k^\times\fleche
\one,
\numeroter \newcount\ffbt \global\ffbt=\formuleno
$$
(cf.~\formule\ratfun\ in the case $F=k(T)$) which can also be derived
more directly by looking at the (tame) symbols corresponding to the
various discrete valuations of $F$.
\th THEOREM 11 (Tate, 1970)
\enonce
For a function field\/ $F$ (in one variable) over a finite field
of\/ $q$ elements, we have\/
$$
\Card\Ker(\lambda)=
(q^2-1)\zeta_F(-1)\Card\Coker(\lambda),
\numeroter
$$
where $\zeta_F$ is the zeta function of\/ $F$.
\endth
\newcount\ffsv \global\ffsv=\formuleno

Let us give a brief sketch of the proof.  
Let $\kbar$ be an algebraic closure of~$k$ and put $\Finfty=F\kbar$ ;
denote the groups of automorphisms of these $\hat{\Z}$-extensions by
$\Gamma=\Gal(\kbar|k)=\Gal(\Finfty|F)$.  Let $D$ be the free
commutative group on the places of $\Finfty$ ({\it i.e.}~on the
rational points of the $\kbar$-curve corresponding to $\Finfty$, or,
equivalently, on the set of discrete valuations of $\Finfty$).  The
maps $\div:\Finftyetoile\fleche D$ ``\thinspace divisor of a
function\thinspace'' and $\deg:D\fleche\Z$ ``\thinspace degree of a
divisor\thinspace'' induce the exact sequences
$$
\one\fleche
\kbaretoile\fleche
\Finftyetoile\fleche
D\fleche
C\fleche
\zero,\quad
\zero\fleche
J(\kbar)\fleche
C\fleche
\Z\fleche
\zero,
$$
of $\Gamma$-modules, where $J$ is the jacobian of $\Finfty$.
Tensoring the first one with $\kbaretoile$ yields the exact sequence
$$
\zero\fleche
\Tor(\kbaretoile,C)\fleche
\kbaretoile\otimes\Finftyetoile\fleche
\kbaretoile\otimes D\fleche
\kbaretoile\fleche
\one
\numeroter \newcount\torseq \global\torseq=\formuleno
$$
($J(\kbar)\otimes\kbaretoile=\zero$ since $J(\kbar)$ is a divisible
and $\kbaretoile$ a torsion group).  On the other hand, taking tame
symbols (cf.~{\it Example\/} 5) at the various discrete valuations of
$\Finfty$ provides an exact sequence of $\Gamma$-modules
$$
\zero\fleche
\Ker(\lambda_\infty)\fleche
\kaytwo(\Finfty)\droite{\lambda_\infty}
\kbaretoile\otimes D\fleche
\kbaretoile\fleche
\one.
\numeroter \newcount\laminf \global\laminf=\formuleno
$$
The map $e:\kbaretoile\otimes\Finftyetoile\fleche\kaytwo(\Finfty)$
which sends $w\otimes f$ to $\symbol(w,f)_\Finfty$ induces a
$\Gamma$-equivariant map of sequences
\formule\torseq$\fleche$\formule\laminf\ which can be shown to be
an isomorphism.
It can be shown further that the map
\formule\ffbt$\fleche$\formule\laminf$^\Gamma$ induced
by the inclusion $F\subset\Finfty$ is also an isomorphism.  Thus
$\Ker(\lambda)$ is isomorphic to $\Ker(\lambda_\infty)^\Gamma$, and
hence to $\Tor(\kbaretoile,C)^\Gamma$, which is the kernel of
$1-\sigma$ acting on $\Tor(\kbaretoile,C)$, where $\sigma\in\Gamma$ is
the (topological) generator $t\mapsto t^q$.  The order of this kernel
is the same as that of the kernel of $1-q\sigma$ acting on~$C$, or,
what comes to the same, $1-q\pi$ acting on the jacobian~$J$ of
$\Finfty$, where $\pi$ is the Frobenius endomorphism.  To conclude the
proof of
\formule\ffsv, it suffices to remark that
$$
\deg(1-q\pi)=
(q^2-1)\zeta_F(-1)(q-1).
$$
 
What is the analogue of \formule\ffsv\ for a number field $F$ ?  The
first difficulty is that the zeta function $\zeta_F$ has a zero of
order $r_2$ ($=$ the number of imaginary places of $F$) at $s=-1$, but
this can be overcome by assuming that $F$ is totally real ($r_2=0$).
Next, one has to interpret the number $q^2-1$ in this context.

So let $F$ be a number field.  Let us restrict ourselves to ultrametric
places, as in \formule\ffbt, and consider only the universal
continuous {\it tame\/} symbols at these places, again as in
\formule\ffbt.  At each place $v$, denote by $k_v$ the residue field ; there
is a canonical surjection $\mu(F_v)\fleche k_v^\times$, which is
an isomorphism for almost all $v$ (the norm residue symbols are tame
at almost all places of $F$).  From \formule\exprec, we deduce an
exact sequence (th.~9)
$$
\kaytwo(F)\droite{\rho}
\sommedir_vk_v^\times\fleche
\one
\numeroter \newcount\nfbt \global\nfbt=\formuleno
$$
which will play the role --- in the totally real case --- played by
\formule\ffbt\  in the function field case.  If we had introduced the
$\kaytwo$ of rings, we could have interpreted $\Ker(\rho)$ as
$\kaytwo(\ogoth)$, where $\ogoth$ is the ring of integers of $F$.

What should be the analogue of the number $q^2-1$ ?  Birch and Tate,
with the help of numerical computations of Atkin, interpreted it as
the largest integer $m$ such that $\Gal(\Fsep|F)$ acts trivially on
$\Z/m\Z(2)$ ; it can also be defined as $\Card\ache^0(F,\Q/\Z(2))$.
They were thus led to the
\th CONJECTURE (Birch-Tate, 1969)
\enonce
For a totally real number field\/ $F$, 
$$
\Card\Ker(\rho)=w_2(F)|\zeta_F(-1)|,
$$
where\/ $w_2(F)$ is the largest\/ $m$ for which\/ $\Gal(\Fsep|F)$ acts
trivially on\/ $\Z/m\Z(2)$.
\endth
%
\goodbreak
The sign of $\zeta_F(-1)$ is $-1$ if and only if $F$ has an odd number
of real places, as follows from the functional equation for $\zeta_F$.
Let us note that Siegel has proved (1969) that $\zeta_F(-1)$ is a
rational number and Serre has proved (1971) that $w_2(F)\zeta_F(-1)$
is an integer.

In the special case $F=\Q$, one has
$\zeta_{\Q}(-1)=-1/12$, $\Card\Ker(\rho)=2$ and $w_2(\Q)=24$.

Tate's result that $K_2(F)\fleche\ache^1(F,\Q/\Z(2))$ is an
isomorphism for a totally real $F$ led Lichtenbaum to formulate his
more general conjectures about $\zeta_F(-m)$, for every odd integer
$m>0$, in terms of the $\Gal(\Fsep|F)$-module $\Q/\Z(m+1)$.



\th THEOREM 12 (Wiles, 1990)
\enonce
The Birch-Tate conjecture is true up to powers of\/~$2$ for every
totally real number field\/ $F$, and indeed exactly true if moreover\/
$F$ is abelian (over\/ $\Q$).
\endth
The proof rests on the main conjecture of Iwasawa theory for all
totally real number fields, which has been proved by Wiles.  However,
for (totally real) {\it abelian\/} number fields, it had been proved
earlier by Mazur \& Wiles.  In this abelian case, the work of
Kolyvagin, Rubin and Thaine, when combined with the earlier work of
Iwasawa, gives a different and simpler proof of the main conjecture.

I thank Prof.~John Coates and R.~Sujatha for their help in getting
these remarks right, and for their critical reading of an earlier
version of this report up to this point.



\section{$\kaytwo$ and special values of $L$-functions
\citer\beilinson(), \citer\bloch(),  \citer\blochgrayson(), \citer\ddjz().}

Let $C$ be a smooth, projective, absolutely connected curve over $\Q$
of genus $g>1$.  Let $N$ be the conductor and $L(C,s)$ the $L$-function
associated to the (continuous, linear) $\Gal(\Qbar|\Q)$-representation
$\ache^1(C_{\Qbar},\Q_l)$ (for any prime number~$l$\/) ; $L(C,s)$ 
converges when the real part of~$s$ is $>{3\over2}$.

\th CONJECTURE  (Hasse-Weil)
\enonce
The function
$$
\Lambda(C,s)={N^{{s\over2}}\over(2\pi)^{gs}}\Gamma(s)^gL(C,s)
$$
admits an analytic continuation to the whole of $\C$ and satisfies
the functional equation $\Lambda(C,s)=w\Lambda(C,s-2)$, with $w=+1$
or\/ $w=-1$.  
\endth

It would follow that $\Lambda(C,0)\in\R^\times$.  The conjecture is
known to be true for modular curves.  It is also true for curves of
genus~$1$, as a result of the seminal work of Wiles and others showing
that all elliptic $\Q$-curves are quotients of (the jacobians of)
modular curves.  We are interested in the special value
$\Lambda(C,0)$.

Let $F=\Q(C)$ be the field of functions of $C$.  It is also the field
of functions of any regular, proper, flat $\Z$-scheme whose generic
fibre is $C$\thinspace; such ``\thinspace integral models\thinspace''
of~$C$ are known to exist.  Fix such a scheme $\Sigma$.  Every
codimension-1 point $P$ of $\Sigma$ gives rise to a discrete valuation
$v_P$ of $F$ and hence ({\it Example\/}~5) to the ``\thinspace tame
symbol\thinspace'' homomorphism $h_P:\kaytwo(F)\rightarrow
k(P)^\times$, where $k(P)$ is the residue field at $P$.  Denote by
$\kaytwo(C,\Z)$ the intersection of the kernels of all these
homomorphisms~; it is independent of the choice of $\Sigma$.  Put
$\kaytwo(C,\Q)=\kaytwo(C,\Z)\otimes_{\Z}\Q$.


\th CONJECTURE (Beilinson)
\enonce
The vector space $\kaytwo(C,\Q)$ is $g$-dimensional.
\endth

Let $X=C(\C)$ be complex analytic curve deduced from $C$.  It is a
compact connected orientable surface of genus $g$ and comes equipped
with a real-analytic involution, induced by complex conjugation.

For $(f,g)\in\Fstar\times\Fstar$ (recall that $F=\Q(C)$), consider the 
(real-analytic) $1$-form
$$
\eta_{f,g}=\log|f|\,d\!\Arg(g)-\log|g|\,d\!\Arg(f)
\numeroter \newcount\iform \global\iform=\formuleno
$$
on the complement in $X$ of the divisors of $f$ and $g$~; it is
bilinear in $f$ and~$g$.  This $1$-form is closed, as $d\eta_{f,g}$
is the imaginary part of $d\log(f)\wedge d\log(g)$, which vanishes.

Let $S$ be a finite subset of $X$ and $\omega$ a closed (smooth)
$1$-form on $X-S$.  For any oriented smooth loop $\gamma$ in
$X-S$, we have the number
$$
(\gamma,\omega)_{X,S}={1\over2\pi}\int_\gamma\omega
$$
which depends only on the class of $\gamma$ in $\ache_1(X-S,\Z)$.  

Let $f,g\in\Fstar$ be such that $f+g=1$ and take $S$ to be complement
in $X$ of the divisors of $f$, $g$.  We have
$(\gamma,\eta_{f,g})_{X,S}=0$, since $\eta_{f,g}=dD(f)$, where $D$ is
the dilogarithm function, a real-analytic function of $z\neq0,1$ in
$\C$ :
$$ 
D(z) = \Arg(1-z) \log |z| - {\rm Im}\left(\int_0^z \log(1-t)
{{dt}\over t} \right). 
$$  

Further, for $s\in S$, let $\gamma_s$ be a smooth loop around $s$ in
$X$, and let $f,g\in\Fstar$.  It can be checked that
$$
(\gamma_s,\eta_{f,g})_{X,S}=\log|t_s(\{f,g\})|,
$$
where $t_s$ is the tame symbol at the place of $F$ determined by $s$.
Thus, if $\{f,g\}$ happens to lie in $\kaytwo(C,\Z)$, then
$(\gamma_s,\eta_{f,g})_{X,S}=0$ for every $s\in S$.  Thus we get a
pairing $\langle\ ,\ \rangle:
\ache_1(X,\Z)\times\kaytwo(C,\Z)\rightarrow\R$. 

The invariants $\ache_1(X,\Z)^+$ under complex conjugation lie in the
left kernel of $\langle\ ,\ \rangle$.  It is the restriction to the
anti-invariants which gives the \bquote regulator map\equote :
$$
\langle\ ,\ \rangle:\ache_1(X,\Q)^-\times\kaytwo(C,\Q)\rightarrow\R,
\numeroter \newcount\regulator \global\regulator=\formuleno
$$
coming from integrating $1$-forms against $1$-cycles, as we have
seen.  Note that $\dim\ache_1(X,\Q)^-=g$.

\th CONJECTURE (Beilinson)
\enonce
The pairing $\langle\ ,\ \rangle$ \formule\regulator\ is perfect.
\endth

Choosing $\Q$-bases for $\ache_1(X,\Q)^-$ and $\kaytwo(C,\Q)$, the
pairing \formule\regulator\ gives a matrix in $\GL_g(\R)$.  The class
of its determinant in $\R^\times\!/\Q^\times$ is independent of the
choice of the $\Q$-bases.

\th CONJECTURE (Beilinson)
\enonce
The determinant of \formule\regulator\ equals\/ $\Lambda(C,0)$ in\/
$\R^\times\!/\Q^\times$. 
\endth

A weak version has been proved for elliptic curves having complex
multiplications.  Similar conjectures have been advanced for curves
over any number field.  Numerical verification has been carried out in
some cases.  

I would like to thank Ramesh Sreekantan and Tim Dokchitser for a
careful reading of this section.

%

\section{Beyond the multiplicative group and the point
\citer\somekawa(), \citer\akhtar()}

$\kaytwo$ is related to numerous other things, and it has been
generalised to numerous other functors.  There is the $\kaytwo$ of a
ring.  There are Milnor's higher $\kay$-groups (for fields), which we
have furtively encountered.  There are Quillen's $\kay$-groups, which
make sense for any scheme, not just for rings and fields.  These
groups can be sheafified, and the cohomology of these sheaves is of
interest.

We have chosen to close by touching upon a recent generalisation due
to Kato, Somekawa and Akhtar.  We shall be even more cryptic here than
elsewhere in this report, as we merely wish to illustrate --- by one
example among many --- how $\kaytwo$ continues to play a central
r{\^o}le in current research.

The idea is to replace the multiplicative group $\Fstar$ of a field
$F$ with an algebraic $F$-group $G$ which is an extension
$$
\zero\rightarrow T\rightarrow G\rightarrow A\rightarrow\zero
$$
of an abelian $F$-variety $A$ by an $F$-torus $T$, i.e.~$G$ is a
semi-abelian variety over $F$.

Let $G_1$, $G_2$, $\ldots$, $G_s$ be semi-abelian $F$-varieties.
Inspired by the presentation of $\kaytwo(F)$ by generators and
relators, the group $\kay(F;G_1,\ldots,G_s)$ is defined by a
presentation a typical generator for which is an element of 
$$
G_1(E)\otimes\cdots\otimes G_s(E)
$$
for a variable finite extension $E$ of $F$ ; they are added as
elements of this group and are subject to two sets of ralators, coming
from the projection formula for $G_i(E)\rightarrow G_i(E')$ for a
morphism $E\rightarrow E'$ of finite extensions, and from a certain
reciprocity law for $G_i(K)\rightarrow\oplus_vG_i(E_v)$ for the
function field $K$ of an $F$-curve, with $v$ running over the places
of $K$ and $E_v$ denoting the residue field at $v$.  When each $G_i$
($1\le i\le s$) is the multiplicative group, one recovers Milnor's
$\kay_s(F)$.

For every $n>1$ invertible in $F$, there is a functorial homomorphism
$c$ which makes the diagram
$$\diagram{
G_1(E)\otimes\cdots\otimes G_s(E)&\hfl{}{}{5mm}&
\ache^1(E,{}_n\!G_1)\otimes\cdots\otimes\ache^1(E,{}_n\!G_s)\cr
\vfl{}{}{5mm}&&\vfl{}{\smile}{5mm}\cr
\kay(F;G_1,\ldots,G_s)\otimes\Z/n\Z&\hfl{c}{}{10mm}&
\ache^s(F,{}_n\!G_1\otimes\cdots\otimes{}_n\!G_s)\cr}
$$
commute for every finite extension $E$ of $F$ ; the map
\formule\bkgalois\ is a particular case.  It is conjectured 
that this $c$ is always injective ; this follows from the Bloch-Kato
conjecture in the case when each $G_i=\Fstar$ ; in this case, the
result is known, as we have seen, if $s=2$ (Merkurjev-Suslin) or if
$n=2^a$ (Voevodsky).

Take $F$ to be a number field, $s=2$, and $G_2=\Fstar$.  There is a
certain local-to-global exact sequence for $\kay(F;G_1,\Fstar)$ of
which Moore's exact sequence \formule\exprec\ is the particular case
$G_1=\Fstar$; the case when $G_1$ is the jacobian of a smooth
projective $F$-curve having a $0$-cycle of degree~$1$ reduces to a
theorem of Bloch and Kato-Saito.

Already Akhtar has further generalised these $\kay$-groups to \bquote
mixed $\kay$-groups\equote.  His construction can be seen as a passage
from $\Spec F$ to a certain number of smooth quasiprojective
$F$-varieties $X_1$, $\ldots$, $X_r$.  A typical generator of the
group $\kay(F;X_1,\ldots,X_r;G_1,\ldots,G_s)$ is an element of
$$
\anought((X_1)_E)\otimes\cdots\otimes\anought((X_r)_E)\otimes
G_1(E)\otimes\cdots\otimes G_s(E)
$$
involving the Chow groups $\anought((X_i)_E)$ (modulo rational
equivalence) of $0$-cycles on $(X_i)_E=X_i\times_FE$, and the groups
of rational points $G_j(E)$, for all finite extensions $E$ of $F$; the
relators come essentially, as before, from the projection formula and
the reciprocity law for local symbols.  The case \hbox{$r=0$} gives
the $\kay$-groups of Kato-Somekawa.  At the other extreme, when
\hbox{$s=0$} and the $X_i$ are projective, one retrieves the group
$\anought(X_1\times_F\cdots\times_F X_r)$ of $0$-cycles modulo
rational equivalence.  The groups $\kay(F;X;G_1,\ldots,G_s)$, where
$X$ is a smooth projective $F$-variety and each $G_i$ is the
multiplicative group $\Fstar$, turn out to be the ``higher Chow
groups'' of $X$ as defined by Bloch.  Interesting computations have
been made over finite fields and local fields, but much remains to be
done.

\bigbreak
\unvbox\bibbox
\bigbreak
{\obeylines\parskip=0pt\parindent=0pt
Chandan Singh Dalawat
Harish-Chandra Research Institute
Chhatnag Road, Jhunsi
{\pc ALLAHABAD} 211\thinspace019, India
\vskip5pt
\tt dalawat@member.ams.org}

\bye

math.AG/0405040

\section{$\kaytwo$ and rational surfaces}
Bloch, Colliot-Thelene, Salberger

\section{$\kaytwo$ and $\Z_l$-extensions}
Coates (Ann. of Math. (2)  95 (1972), 99--116; MR 50 #12971)

MR1144791 (92m:11121)
Kersten, Ina(D-BLF)
$K\sb 2$ und $ Z\sb p$-Erweiterungen von $\bold Q(\zeta\sb {p\sp
r})$. (German. English summary) [$K\sb 2$ and $\bold Z\sb
p$-extensions of $\bold Q(\zeta\sb {p\sp r})$]
Mathematische Wissenschaften gestern und heute. 300 Jahre
Mathematische Gesellschaft in Hamburg, Teil 2.
Mitt. Math. Ges. Hamburg 12 (1991), no. 2, 347--362.

Consider cyclotomic fields $K=\bold Q(\zeta_{p^r})$ of $p$-power roots
of unity, $p$ an odd prime. Let $\tilde K$ be the compositum of all
$\bold Z_p$-extensions of $K$, $A_K=\{a\in K^*$: $K(\root
p\of{a})\subset\tilde K\}$, $C_K=\{c\in K^*$: $\{c,\zeta_p\}=1$ in
$K_2(K)\}$. Concerning a question of Coates, it was proved by
R. Greenberg \ref[Amer. J. Math. 100 (1978), no. 6, 1235--1245; MR
80d:12002] that $A_K=C_K$ for $r$ sufficiently large. In this paper it
is proved that $A_K/K^{*p}=U/U^p=C_K/K^{*p}$ for $r$ sufficiently
large, where $U$ is the group of $p$-units of $K$.

MR1176424 (93f:11089)
Kolster, Manfred(3-MMAS-MS)
$K\sb 2$ of rings of algebraic integers.
J. Number Theory 42 (1992), no. 1, 103--122.
11R70 (19C99 19F99)

Let $F$ be an algebraic number field with ring of integers ${\scr
O}_F$. For odd prime numbers $p$ and totally real $F$, J. Coates has
related \ref[Ann. of Math. (2) 95 (1972), 99--116; MR 50 #12971] the
$p$-primary part of the tame kernel $K_2({\scr O}_F)$ to fixed
elements of a Tate twisted ideal class group that occurs naturally in
Iwasawa theory. In the first part of this paper, the author
generalizes Coates' result to arbitrary number fields.

\section{$kaytwo$ and congurence subgroups}
cf. Milnor's book, p. xi

\section{$kaytwo$ and isotopies}
cf. Milnor's book, p. x

MR0368036 (51 #4278) Hatcher, Allen E.  Pseudo-isotopy and $K\sb{2}$.
Algebraic $K$-theory, II: "Classical" algebraic $K$-theory and
connections with arithmetic (Proc. Conf., Seattle Res. Center,
Battelle Memorial Inst., 1972), pp. 328--336. Lecture Notes in Math.,
Vol. 342, Springer, Berlin, 1973.

MR0353337 (50 #5821)
Hatcher, Allen; Wagoner, John
Pseudo-isotopies of compact manifolds.
With English and French prefaces. Astérisque, No. 6.
Société Mathématique de France, Paris, 1973. i+275 pp.

\section{$kaytwo$ and affine groups}
Brylinski, Deligne, Prasad, Raghunathan

MR1441006 (98d:14056)
Deligne, P.(1-IASP)
Extensions centrales de groupes algébriques simplement connexes et
cohomologie galoisienne. (French) [Central extensions of simply
connected algebraic groups and Galois cohomology]
Inst. Hautes Études Sci. Publ. Math. No. 84, (1996), 35--89 (1997).
14L15 (14F20 19C09)

MR1896177 (Review)
Brylinski, Jean-Luc(1-PAS); Deligne, Pierre(1-IASP)
Central extensions of reductive groups by $\bold K\sb 2$.
Publ. Math. Inst. Hautes Études Sci. No. 94, (2001), 5--85.

\section{$kaytwo$ and galoisian descent}
MR1220427 (94i:11094)
Kahn, Bruno(F-PARIS7)
Descente galoisienne et $K\sb 2$ des corps de
nombres. (French. English, French summary) [Galois descent and $K\sb
2$ of number fields]
$K$-Theory 7 (1993), no. 1, 55--100.

Let $E/F$ be a finite Galois extension with group $G$ and let
$f_{E/F}\colon K_2(F)\to K_2(E)^G$ denote the canonical map. The main
results are natural isomorphisms of the form ${\rm Ker}(f_{E/F})\cong
H^1(G;K^{\rm ind}_3 (E_0))$ and ${\rm Coker}(f_{E/F})\cong
H^2(G;K^{\rm ind}_3 (E_0))$, where $E_0$ denotes the field of
constants and $K^{\rm ind}_3(E)$ denotes the indecomposable
$K_3$-group of $E$. The proof uses the hypercohomology of the $\Gamma
(2)$-complex of S. Lichtenbaum \ref[Invent. Math.  88 (1987), no. 1,
183--215; MR 88d:14011]. Several applications are given to function
fields of varieties and to number fields. In the latter situation
several of the author's results were obtained independently by
J. Brinkhuis \ref[in Algebraic $K$-theory, number theory, geometry and
analysis (Bielefeld, 1982), 13--28, Lecture Notes in Math., 1046,
Springer, Berlin, 1984; MR 85j:12008] and M. Kolster \ref[J. Pure
Appl. Algebra 74 (1991), no. 3, 257--273; MR 93c:19004].

\section{$\kaytwo$ and intersection theory}
Bloch, Quillen, Grayson

Information about the Chow ring.

Chow group of codimension~2 isomorphic to H^2 of the sheafified K_2, etc.

\section{$\kaytow$ and the dilogarithm}

MR0578856 (82f:14009) Bloch, Spencer Applications of the dilogarithm
function in algebraic $K$-theory and algebraic geometry.  Proceedings
of the International Symposium on Algebraic Geometry (Kyoto Univ.,
Kyoto, 1977), pp. 103--114, Kinokuniya Book Store, Tokyo, 1978.  

The author starts by tensoring the exponential sequence with $ C^*$:
$$ 0\fleche\{ C^*/\text{torsion}\}\fleche C\otimes_{ Z}
C^*\fleche C^*\otimes_{ Z} C^*\fleche 0. $$ Now observe that
by results of Matsumoto we have an exact sequence: $$ 0\fleche B(
C)\fleche A( C)\overset\lambda\to\fleche C^*\otimes_{ Z}
C^*\fleche K_2( C)\fleche 0 $$ where $A( C)$ is the free
abelian group on the set $ C^*-\{1\}$ and if $a\in C^*-\{1\}$,
$\lambda[a]=a\otimes(1-a)$. To obtain a sequence involving $K_2(C)$
rather than $ C^*\otimes_Z C^*$ the author lifts $\lambda$ to a map
$\varepsilon\colon A( C)\fleche C\otimes_{ Z} C^*$ defined by: $$
\varepsilon[a]=\log(1-a)\otimes a+2\pi
i\otimes\exp\left(\frac{-1}{2\pi
i}\int_0^a\log(1-t)t\inv \,dt\right). $$ $v( C)$ is defined as $
C\otimes C^*/\varepsilon(A( C))$. Finally, the author sketches a proof
of the remarkable fact that $\tilde\Delta^*=\varepsilon(B( C))$ is a
(nontrivial) subgroup of $\overline{ Q}^*/\text{torsion}$.

\section{$\kaytwo$ and algebraic cycles}

Nesterenko-Suslin, Totaro

MR1187705 (94d:19009) Totaro, Burt(1-CHI) Milnor $K$-theory is the
simplest part of algebraic $K$-theory. (English. English summary)
$K$-Theory 6 (1992), no. 2, 177--189.  

Let $F$ be a field and denote by ${\rm CH}^p(F,n)$ S. J. Bloch's
higher Chow groups of $F$ \ref[Adv. in Math.  61 (1986), no. 3,
267--304; MR 88f:18010]. These groups were computed by
Yu. P. Nesterenko and A. A. Suslin for $p\geq n$ \ref[Izv. Akad. Nauk
SSSR Ser. Mat.  53 (1989), no. 1, 121--146; MR 90a:20092]: they are
zero for $p>n$, and for $p=n$ it is the Milnor $K_n$-group of $F$. In
the paper under review, the author gives another proof of this
result.  It is quite explicit and very nice.

As preparation for the proof, the author first replaces Bloch's
complex with another complex constructed using cubes instead of
simplices. This makes it easier to describe the product structure on
the higher Chow groups. Then he replaces affine $n$-space with $(
P^1-\{1\})^n$. This makes the computation of ${\rm CH}^1(F,1)$
easier. He defines a map $F^*\to {\rm CH}^1(F,1)$ in a very explicit
way and then uses the product structure to extend this to a map:
$K^M_i(F)\to {\rm CH}^i(F,i)$. Explicit curves are produced to help
verify the necessary relations in order to prove the result. The proof
of the main theorem is different for infinite fields than for finite
fields.

\section{$\kaytwo$ and etale cohomology}
eg. MR0553999 (81i:12016)
Soulé, C.
$K$-théorie des anneaux d'entiers de corps de nombres et cohomologie
étale. (French)
Invent. Math. 55 (1979), no. 3, 251--295.

\section{$\kaytwo$ and cristalline cohomology}
eg. MR0488288 (81j:14011)
Bloch, Spencer
Algebraic $K$-theory and crystalline cohomology.
Inst. Hautes Études sci. Publ. math. No. 47, (1977), 187--268 (1978).

\section{$\kaytwo$ and motivic cohomology}
Beilinson, Lichtenbaum, Voevodsky

MR1807268 (2003a:14031) Geisser, Thomas(J-TOKYO); Levine, Marc(1-NORE)
The Bloch-Kato conjecture and a theorem of
Suslin-Voevodsky. (English. English summary) J. Reine Angew. Math. 530
(2001), 55--103.

There is a list of conjectural properties of motivic cohomology
\ref[see S. Lichtenbaum, in Motives (Seattle, WA, 1991), 303--313,
Proc. Sympos. Pure Math., 55, Part 1, Amer. Math. Soc., Providence,
RI, 1994; MR 95a:19007], several of which have already been
established. For instance, it is known that $H^p_M(F,{\Bbb
Z}(p))\cong K^M_p(F)$ for fields $F$ where $K^M_p$ is the $p$th
Milnor $K$-theory group.

\section{$\kaytwo$ and regulators}

MR1809627 (2002c:14036)
Besser, Amnon(4-DRHM)
Syntomic regulators and $p$-adic integration. II. $K\sb 2$ of
curves. (English. English summary)
Proceedings of the Conference on $p$-adic Aspects of the Theory of
Automorphic Representations (Jerusalem, 1998).
Israel J. Math. 120 (2000), , part B, 335--359.
 
MR0562648 (81h:12011)
Bloch, S.
Algebraic $K$-theory and zeta functions of elliptic curves.
Proceedings of the International Congress of Mathematicians (Helsinki,
1978), pp. 511--515,
Acad. Sci. Fennica, Helsinki, 1980.
12A62 (14K07)

MR0862631 (88f:11061)
Bloch, S.(1-CHI); Grayson, D.(1-IL)
$K\sb 2$ and $L$-functions of elliptic curves: computer calculations.
Applications of algebraic $K$-theory to algebraic geometry and number
theory, Part I, II (Boulder, Colo., 1983), 79--88,
Contemp. Math., 55,
Amer. Math. Soc., Providence, RI, 1986.

The authors describe computer calculations which were used to test a
conjecture of Bloch \ref[ Proceedings of the International Congress of
Mathematicians (Helsinki, 1978), 511--515, Acad. Sci. Fennica,
Helsinki, 1980; MR 81h:12011] and \n A. Be\u\i linson\en
\ref[Funktsional. Anal. i Prilozhen. 14 (1980), no. 2, 46--47; MR
81k:14020] relating the regulator on $K_2$ of a modular elliptic curve
$E_{/\bold Q}$ to the value at $s=2$ of the $L$-function of $E$. When
$E$ has complex multiplication a proof of the conjecture has been
given by \n D. Rohrlich\en \ref[ Number theory (Montreal, Que., 1985),
371--387, Amer. Math. Soc., Providence, R.I., 1987]. However, in case
there are primes of multiplicative reduction, the computer
calculations of this paper suggest that the conjectures need to be
modified by replacing $K_2(E_{\bold Q})$ by $K_2(E_{\bold Z})$ where
$E_{\bold Z}$ is the Néron model of $E$.

MR0575206 (81k:14020)
Be\u\i linson, A. A.
Higher regulators and values of $L$-functions of curves. (Russian)
Funktsional. Anal. i Prilozhen. 14 (1980), no. 2, 46--47.
14G10 (12A70)

Let $X$ be a nonsingular algebraic curve over the field of complex
numbers $ C$. The author defines a mapping from $K_2(X)$ to $H^1(X,
C^*)$, which generalizes the tame symbol (the same construction was
obtained independently by P. Deligne, but was not published). Then a
regulator is determined for elliptic and modular curves in terms of
the value of their $L$-functions in points 1 and 2. The author gives a
construction of $ C^*$-valued higher regulators for points and
formulates a rigidity theorem for Borel's regulators. These results
are generalizations of some results of S. Bloch (preprints).

The author conjectures that if $X$ is defined over the field of
rational numbers $ Q$, then the $L$-functions at the point 2 can be
expressed up to a rational factor by logarithms of the moduli of
periodic elements of $K_2$ defined over $ Q$.

MR0948100 (89k:11041)
Coleman, Robert(1-CA); de Shalit, Ehud(IL-HEBR)
$p$-adic regulators on curves and special values of $p$-adic $L$-functions.
Invent. Math. 93 (1988), no. 2, 239--266.

The purpose of this paper is to define for smooth complete curves
$C/\overline{\bold Q}_p$ a $p$-adic regulator map on $K_2(C)$ and, in
case $C$ is an elliptic curve with complex multiplication, to relate
the map to a special value of the $p$-adic $L$-series of $C$. More
specifically, following earlier work of Coleman \ref[same journal 69
(1982), no. 2, 171--208; MR 84a:12021], the authors develop a theory
of $p$-adic integration on $C$. Using this, they define a map $r\:K_2
(\overline{\bold
Q}_p(C))\to\roman{Hom}(H^0(C,\Omega^1_C),\overline{\bold Q}_p)$, whose
value on the Steinberg symbol $\{f,g\}$ is the functional
$\omega\mapsto\int_{(f)} Log(g)·\omega$, the integration being over
the divisor of $f$ and $ Log$ a fixed branch of the $p$-adic
logarithm.

Let $C$ be an elliptic curve defined over $\bold Q$ with complex
multiplication by the imaginary quadratic field $K$, and let $\psi$ be
the associated Hecke character of $K$. Suppose $p=P\overline P$ splits
in $K$ and $C$ has good ordinary reduction at $p$. If $f,g$ are
functions on $C$ whose divisors are supported on torsion points, it is
shown that $(1-1/p\psi(P))r(f,g)(\omega)\in\Omega_p(\omega)L_p(\psi)
\bold Q$, where $\Omega_p(\omega)$ is the $p$-adic period of $\omega$
and $L_p$ is the $p$-adic $L$-function of Manin-Vishik and Katz. This
is the $p$-adic version of a result of \n S. Bloch\en \ref[in
Proceedings of the International Congress of Mathematicians (Helsinki,
1978), 511--515, Acad. Sci. Fennica, Helsinki, 1980; MR 81h:12011; see
also D. Rohrlich, in Number theory (Montréal, PQ, 1985), 371--387,
Amer. Math. Soc., Providence, RI, 1987; MR 89b:11054].

\section{$K_2$ and regular local rings}
math.KT/0402190
On K_2 of 1-dimensional local rings
Authors: Amalendu Krishna

\begin{conjecture} ({\bf Geller})
Let $A$ be a local one-dimensional domain with field of fractions $F$. Then
$A$ is regular if and only if the map $K_2(A) \lra K_2(F)$ is injective.
\end{conjecture}
! Maybe this $K_2$ is the Quillen group and not the Milnor group.

This conjecture was verified by Geller ([G]) when $A$ is noetherian,
equicharacteristic, characteristic zero, and is also seminormal with
finite normalisation. In the same article, Dennis and Sherman verify
it for cuspidal rings of the type $K[t^2, t^3]$ as described
above. The conjecture is still unknown in almost all other cases. Our
first aim in this paper is to formulate an Artinian version of this
conjecture, and to show that this conjecture implies Geller's
conjecture.

\bibitem[G]{ }  S. Geller,  {\it A note on injectivity of lower $K$-groups
for integral domains\/}, with appendix by R. Dennis and C. Sherman,
Contemporary Math., {\bf 55}, Vol. II, (1986), 437-447.
MR0862647

\section{$\kaytwo$ and Galois modules}

MR1641555 (99g:11140)
Chinburg, T.(1-PA); Kolster, M.(3-MMAS-MS); Pappas, G.(1-PRIN);
Snaith, V.(4-SHMP) 
Galois structure of $K$-groups of rings of integers. (English. English summary)
$K$-Theory 14 (1998), no. 4, 319--369.

The analytic class number formula for a number field $N$ expresses the
product of the class number and the regulator of $K$ as the residue of
the Dedekind zeta function of $N$, times some explicit and easily
computable factors. (The regulator "measures" in some way the unit
group of the ring of integers in $N$.) This classical link between
class numbers and unit groups (via regulators) has been the subject of
far-reaching generalizations, both established and conjectural, and
this paper is an important contribution. The objects of study are the
$K$-groups $K_{2n+1}(O_N)$ (odd-numbered) and $K_{2n}(O_N)$
(even-numbered) of the ring of integers $O_N$, where $N/K$ is now a
$G$-Galois extension of number fields (or of function fields, see
Chapter 7). It is helpful to think of the odd-numbered $K$-groups as
"twisted versions" of the unit group, and of the even-numbered ones as
"twisted versions" of the class group. Indeed, $K_1(O_N)$ is precisely
the unit group, and $K_0(O_N)$ is isomorphic to ${\rm Cl}_N\oplus Z$.
All these $K$-groups are considered as $G$-modules. The difference
$D=(K_{2n+1}(O_N)-{}(K_{2n}(O_N)))$ should be something simple, viewed
in a suitable Grothendieck group.

One should say what happens for (K_2)-(K_1).

\section{misc.}

MR0335522 (49 #303)
Krusemeyer, Mark I.
Fundamental groups, algebraic $K$-theory, and a problem of Abhyankar.
Invent. Math. 19 (1973), 15--47.

L'auteur présente une application intéressante de la $K$-théorie
algébrique en géométrie algébrique, où il utilise le foncteur $K_2$ de
Milnor et le groupe fondamental $\pi_1(\text{SL})$ du foncteur SL pour
résoudre négativement un problème de S. S. Abhyankar [ Algebraic space
curves, Presses Univ. Montréal, Montreal, Que., 1971]. Ce problème
peut se réduire au problème suivant que l'auteur appelle
$(\alpha,\beta,\gamma)$-problème: Soient $\alpha,\beta,\gamma\in
k[X,Y]$ tels que $(\alpha,\beta,\gamma)k[X,Y]$ soit un idéal à deux
générateurs, $k$ étant un corps commutatif algébriquement
clos. Existe-t-il deux éléments $a,b\in k[X,Y]$ tels que l'on ait
$(\alpha,\beta,\gamma)k[X,Y]=(\alpha+a\gamma,\beta+b\gamma)k[X,Y]$?

MR0371891 (51 #8108)
Bloch, Spencer
$K\sb{2}$ of Artinian $Q$-algebras, with application to algebraic cycles.
Comm. Algebra 3 (1975), 405--428.
14C15 (14C25 18F25)

MR1753108 (2002b:14007)
Raskind, Wayne(1-SCA); Spiess, Michael(D-HDBG)
Milnor $K$-groups and zero-cycles on products of curves over $p$-adic
fields. (English. English summary)
Compositio Math. 121 (2000), no. 1, 1--33.

DOCUMENTA MATHEMATICA, Extra Volume: Kazuya Kato's Fiftieth Birthday
(2003), 387-442 Takako Fukaya

Coleman Power Series for $K_2$ and $p$-Adic Zeta Functions of Modular Forms

For a usual local field of mixed characteristic $(0,p)$, we have the
theory of Coleman power series \cite{Co}. By applying this theory to
the norm compatible system of cyclotomic elements, we obtain the
$p$-adic Riemann zeta function of Kubota-Leopoldt \cite{KL}. This
application is very important in cyclotomic Iwasawa theory.\par In
\cite{Fu}, the author defined and studied Coleman power series for
$K_2$ for certain class of local fields. The aim of this paper is
following the analogy with the above classical case, to obtain
$p$-adic zeta functions of various cusp forms (both in one variable
attached to cusp forms, and in two variables attached to ordinary
families of cusp forms) by Amice-Vélu, Vishik, Greenberg-Stevens, and
Kitagawa,... by applying the $K_2$ Coleman power series to the norm
compatible system of Beilinson elements defined by Kato \cite{Ka7} in
the projective limit of $K_2$ of modular curves.

MR1771574 (2001k:14020)
Esnault, Hélène(D-ESSN); Levine, Marc(1-NORE)
The Steinberg curve. (English. English summary)
Amer. J. Math. 122 (2000), no. 4, 783--804.
14C25 (14C15 14F43 19C20 19E15)

The basic idea of the proof is astonishingly simple. Consider
surjective analytic maps $p\colon \, C^{\times}\to E$, $p'\colon \,
C^{\times}\to E'$ coming from analytic uniformizations of the elliptic
curves in question. Then for any $u,v\in C^{\times}$ the formula
$p(u)*p'(v)\coloneq (p(u),p'(v))-(p(u),0)-(0,p'(v))-(0,0)$ defines a
degree 0 element in ${\rm CH}_0(X)=$\break ${\rm CH}^2(X)$. One
verifies that for any $u\neq 0$ the "Steinberg cycles"\break
$p(u)*p'(1-u)$ lie in the kernel of the cycle map $\phi$.

MR1629866 (99i:19001)
Kurihara, Masato(J-TOKYM)
The exponential homomorphisms for the Milnor $K$-groups and an
explicit reciprocity law. (English. English summary)
J. Reine Angew. Math. 498 (1998), 201--221.
19D45 (11S70)

Let $K$ be a discrete valuation field of characteristic $0$ with ring
of integers $O_K$ and residue field $F$ of characteristic $p>0$. Let
$v_K$ denote the normalized additive valuation on $K$ and $e=v_K(p)$
the absolute ramification index. Generalizing the classical
exponential homomorphism $\exp\colon O_K\to O_K^*$ the author
constructs for each $\eta \in O_K$ such that $v_K(\eta) \geq 2e/(p-1)$
and all $q>0$ exponential homomorphisms $\exp_{\eta,q}\colon
\Omega_A^{q-1} \fleche K_q^M(A)^{\wedge}$ such that
$$a\frac{db_1}{b_1} \wedge \cdots \wedge \frac{db_{q-1}}{b_{q-1}}
\mapsto \{\exp(\eta \,a),b_1,\cdots,b_{q-1}\}$$ for $a \in A$, $b_i
\in A^*$ provided that $F \not= \bold F_2$. Here $A$ is equal either
to $O_K$ or to a formal power series ring over $O_K$, $\Omega_A^{q-1}
$ is the\break $(q-1)$st exterior power of the module of absolute
K\"{a}hler differentials and $$K_q^M(A)^{\wedge} =
\varprojlim\,K_q^M(A)/p^n\,K_q^M(A)$$ is the $p$-adic completion of
the Milnor $K$-group $K_q^M(A)$ of $A$. In the special case in which
$p=2$ and $F=\bold F_2$ a similar construction yields exponential
homomorphisms with values in the $2$-adic completion of the Milnor
$K$-theory of the ring $A[\frac{1}{2}]$.

As an application of the higher exponential homomorphisms the author
gives a simple proof of the explicit reciprocity law for local fields
of S. Sen \ref[J. Reine Angew. Math. 313 (1980), 1--26; MR 81a:12018]
and generalizes the explicit reciprocity laws to higher-dimensional
local fields.

MR1136845 (92k:11068)
Kato, Kazuya(J-TOKYOS)
The explicit reciprocity law and the cohomology of
Fontaine-Messing. (English. French summary)
Bull. Soc. Math. France 119 (1991), no. 4, 397--441.
11G45 (11R37 14F20 14F30 19D45)

The author provides a factorization over crystalline cohomology of a
symbol map, for certain types of regular base-rings $A$. The symbol
goes from Milnor $K$-theory to étale cohomology. The map from
crystalline to étale cohomology is given by a simple direct
definition, and is apparently equal to the map defined by Fontaine and
Messing (who have announced such a result, but have not published a
proof). The author is very careful to give complete proofs or
references for all the results. Unfortunately, the long delay in
publishing the paper has not been used to update its contents, or at
least all parts of the bibliography.

MR1756261 (2001i:11139)
Nakamura, Jinya(J-TOKYOGM)
On the structures of the Milnor $K$-groups of some complete discrete
valuation fields. (English. English summary)
$K$-Theory 19 (2000), no. 3, 269--309.
11S70 (19D45)

Let $L$ be a complete discrete valuation field with valuation ring
$O_F$ and valuation ideal $\germ m$. The unit group $L^*$ of $L$ has a
natural filtration $\{U_L^i\}_{i \geq 0}$, where $U_L^0 = O_L^*$ and
$U_L^i = 1 + \germ m^i_L$ for $i \geq 1$. The graded quotients
$U_L^i/U_L^{i+1}$ are known to be isomorphic to the unit group of the
residue field if $i = 0$, and to the additive group of the residue
field otherwise.

An analogous filtration can be defined on the $q$-th Milnor $K$-group
$K_q^M(L)$ ($q \geq 1$): For $i \geq 0$ let $U_q^i(L)$ denote the
subgroup of $K_q^M(L)$ generated by all symbols
$\{x_1,x_2,\cdots,x_q\}$ with $x_1 \in U_L^i$ and $x_i \in L^*, i \geq
2$. The structure of the graded quotients ${\rm
gr}_q^i(L)=U_q^i(L)/U_q^{i+1}(L)$ is known in the case of equal
characteristic due to results of J. Graham \ref[in Algebraic
$K$-theory, II: ``Classical'' algebraic $K$-theory and connections
with arithmetic (Proc. Conf., Battelle Memorial Inst., Seattle, Wash.,
1972), 474--486. Lecture Notes in Math., 342, Springer, Berlin, 1973;
MR 51 #442] and S. Bloch\ \ref[Inst. Hautes Études
Sci. Publ. math. No. 47, (1977), 187--268 (1978); MR 81j:14011].

In the case of mixed characteristic $(0,p)$ much less is
known. M. Kurihara\ \ref[Invent. Math. 93 (1988), no. 2, 451--480; MR
89k:11118] determined the graded quotients ${\rm gr}_q^i(L)$ in the
case that $L$ is absolutely ramified, but in general the graded
quotients are only known in the range $0 \leq i \leq ep/(p-1)$, where
$e = {\rm ord}_L(p)$ is the absolute ramification index
\ref[cf. S. J. Bloch\ and K. Kat\=o, Inst. Hautes Études
Sci. Publ. Math. No. 63 (1986), 107--152; MR 87k:14018]. The main
result of the paper under consideration gives a complete calculation
of all graded quotients ${\rm gr}_q^i(L)$ for certain fields $L$ which
are not absolutely unramified. Typical examples arise as follows: Let
$p$ be an odd prime, and let $R$ be the completion of the localization
of the polynomial ring $\bold Z_p[T_1,\cdots,T_n]$ at the prime ideal
$(p)$. Let $K$ be the fraction field of $R$ and take $L =
K(\zeta_p,\root \uproot 2 p \of {1-\zeta_p})$ or $L = K(\zeta_p, \root
\uproot 2 p \of {-p})$.

The computations of the graded quotients ${\rm gr}_q^i(L)$ for the
fields $L$ considered are in terms of exact sequences involving
exterior powers of the module of absolute differentials of the residue
field.

MR0369322 (51 #5557)
Coates, J.; Sinnott, W.
An analogue of Stickelberger's theorem for the higher $K$-groups.
Invent. Math. 24 (1974), 149--161.
12A70

Let $F$ be an abelian field with ring of integers ${\germ o}_F$. Let
$f$ denote the conductor of $F$ and $G={\rm Gal}(F/ Q)$. As a
generalization of the classical Stickelberger elements, J. Coates and
W. Sinnott defined twisted Stickelberger elements \ref[Invent. Math.
24 (1974), 149--161; MR 51 #5557] in the following manner: Let $j$ be
an integer $\geq 0$. For each integer $b\geq 1$ such that $(b,f)=1$,
define a twisted Stickelberger element $$S_j(b)=w_{j+1}(
Q)(b^{j+1}-(b,F))·\sum\Sb a\bmod f\\
(a,f)=1\endSb\zeta_F(a,-j)(a,F)^{-1}.$$ Here $(\;,F)$ denotes the
Artin symbol, $\zeta_F$ the partial $\zeta$ function and $w_{j+1}( Q)$
is the order of the Galois cohomology group $H^0(F, Q_p/ Z_p(j+1))$.

The $S_j(b)$ are in fact in $ Z[G]$ and therefore generate an ideal
$S_j(F)$ in $ Z[G]$. Coates and Sinnott conjectured that $S_j(F)$
annihilates the Quillen $K$-group $K_{2j}({\germ o}_F)$, proving it
for $j=1$ up to 2-torsion. (Kolster, review of MR1154596).

MR1713135 (2001i:12004)
Efrat, Ido(IL-BGUN-CS); Fesenko, Ivan(4-NOTT)
Fields Galois-equivalent to a local field of positive characteristic.
Math. Res. Lett. 6 (1999), no. 3-4, 345--356.
12F12 (12J10)

Let $F= F_q ((T))$ be the Laurent series field with coefficients in
the finite field $ F_q$. The structure of the absolute Galois group
$G$ of $F$ has been determined by H. Koch \ref[Math. Nachr.  35
(1967), 323--327; MR 37 #5189]. The problem considered here is: What
can one say of a field $K$ whose absolute Galois group is isomorphic
to $G$? The authors prove that $K$ possesses a Henselian valuation $v$
with residue field of characteristic $p$ and satisfying other
properties. The ingredients of the proof are powerful methods of field
theory and in particular the construction of valuations from subgroups
of the multiplicative group and Milnor $K_2$. The authors also give
examples of such fields $F$ whose valuations have various properties.

MR1713134 (2001i:12011)
Efrat, Ido(IL-BGUN-CS)
Construction of valuations from $K$-theory. (English. English summary)
Math. Res. Lett. 6 (1999), no. 3-4, 335--343.
12J25 (19C20)

One of the objectives of field theory is to obtain information
(characterization) of a field $K$ from the its abolute Galois
group. Let $p$ be a prime number. One of the methods to do that is to
construct valuations from subgroups $T$ of the quotient of the
multiplicative group $K^*$ by its $p$-powers. The subgroup $T$ has to
satisfy appropriate properties with respect to Milnor $K_2$. In this
expository paper, the author gives an elementary approach to this
construction.

MR1936475 (2003k:14042)
Pablos Romo, Fernando(E-SALA)
On the tame symbol of an algebraic curve. (English. English summary)
Comm. Algebra 30 (2002), no. 9, 4349--4368.
14H99 (19F15)

MR1809368 (2002g:12003)
Efrat, Ido(IL-BGUN-CS)
The local correspondence over absolute fields: an algebraic approach.
Internat. Math. Res. Notices 2000, no. 23, 1213--1223.
12F10

A field $K$ is called absolute if it is finitely generated over its
prime field. The zero-dimensional case of Grothendieck's "anabelian
conjecture" says that an infinite absolute field is essentially
determined by its absolute Galois group $G_K$. Finally this was proven
by Pop using model-theoretic methods.

The note under review gives a purely algebraic proof of the needed
local correspondence. This new proof is based on an elementary
construction of valuations on fields with "large" second Milnor
$K$-group.

Reviewed by Martin Epkenhans

From Romyar Sharifi's homepage :

A cup product pairing: $R_K = Z[\mu_p, 1/p]$ is the ring of
$p$-integers of $K=\Q(\mu_p)$. $E_K = R_K^\times$ is the group of
$p$-units of $K$. McCallum and I defined a pairing $( , )_K :
E_K\times E_K\rightrrow A_K\otimes\mu_p$ (where $A_K$ is the $p$-part
of the class group of $K$), which arises from the cup product in
\´etale (or Galois) cohomology $H^1(\Spec R_K,
\mu_p)^{\otimes2}\rightarrow H^2(\Spec R_K,\mu_p^{\otimes2})$.

Conjecture (McCallum-S). $( , )_K$ is surjective. 

Theorem (S). $( , )_K$ is surjective for $p < 1000$.

Milnor K-groups: Define $K^M_2(R_K) = {E_K\otimes EK \over
x\otimes(1-x) | x, 1-x \in E_K}$. 

We have a canonical homomorphism $K^M_2(R_K)\rightarrow K_2(R_K)$,
where $K_2(R_K)$ is the usual algebraic $K_2$-group. 

Remark. If $R_K$ is replaced by any field and $E_K$ by its
multiplicative group, the above map is an isomorphism
(Matsumoto). 

Surjectivity of $( , )_K$ can be reinterpreted as the
following equivalent statement.

Conjecture (McCallum-S). The map $K^M_2\otimes\Z_p\rightarrow
K_2(R_K)\otimes\Z_p$ is surjective. 

\\
Paper: math.KT/0511211
Date: Tue, 8 Nov 2005 17:46:30 GMT   (10kb)

Title: A criterion for cohomological dimension
Authors: Karen Acquista
Comments: 9 pages
Subj-class: K-Theory and Homology; Number Theory
\\
 We give a criterion for the cohomological dimension of a field, involving
norm maps on Milnor K-theory; this criterion was originally formulated by Kato.
The theorem we prove is a generalization of a theorem in Serre's book on Galois
cohomology.
\\ ( http://arXiv.org/abs/math/0511211 ,  10kb)

1. math.AG/0511190 [abs, ps, pdf, other] :
    Title: K_2 of elliptic surface minus singular fibers and
    q-expansion of Beilinson's regulator 
    Authors: Masanori Asakura
    Comments: 49 pages
    Subj-class: Algebraic Geometry; Number Theory
    MSC-class: 14J27; 19D45; 19F27
     
2. math.NT/0502222 [abs, ps, pdf, other] :
    Title: Surjectivity of $p$-adic regulator on $K_2$ of Tate curves
    Authors: Masanori Asakura
    Subj-class: Number Theory; Algebraic Geometry
    MSC-class: 19F27

Greenberg, Peter
Area preserving pl homeomorphisms and relations in $K_2$. Annales de
l'institut Fourier, 48 no. 1 (1998), p. 133-148 
Full text djvu | pdf | Reviews MR 99d:19001 | Zbl 0904.19001

Abstract 

To any compactly supported, area preserving, piecewise linear
homeomorphism of the plane is associated a relation in of the smallest
field whose elements are needed to write the homeomorphism.

Using a formula of J. Morita, we show how to calculate the relation,
in some simple cases. As applications, a ``reciprocity" formula for a
pair of triangles in the plane, and some explicit elements of torsion
in of certain function fields are found.

MR1614890 (99d:19001)
Greenberg, Peter
Area preserving ${\rm pl}$ homeomorphisms and relations in $K\sb
2$. (English. English, French summary) 
Ann. Inst. Fourier (Grenoble) 48 (1998), no. 1, 133--148.
19C20 (52B45)
 
Let ${\rm SPL}_C\bold R^2$ denote the group of compactly supported,
area-preserving PL-homeomorphisms of the plane $\bold R^2$. Using a
construction from the author and V. Sergiescu \ref[$K$-Theory 9
(1995), no. 6, 529--544; MR1361581 (97c:57037)] the author constructs
a function, $D_g\colon\bold R^2\to K_2(\bold R)$, for each $g\in{\rm
SPL}_2\bold R^2$ such that (i) $D_g(v)\neq 1$ for only finitely many
$v\in\bold R^2$, (ii) $D_{gh}(v)=D_g(hv)D_h(v)$ and $\prod_{v\in\bold
R^2}D_g(v)=1$.

Using this result, geometrical constructions are used to construct
elements and relations in $K_2$ of $\bold R$ and of function
fields. The methods are reminiscent of the scissors congruence methods
of S. Lichtenbaum \ref[in Algebraic $K$-theory and algebraic number
theory (Honolulu, HI, 1987), 151--157, Amer. Math. Soc., Providence,
RI, 1989; MR0991980 (90e:20030)] and the polylogarithmic results of
A. B. Goncharov \ref[Adv. Math. 114 (1995), no. 2, 197--318; MR1348706
(96g:19005)].

Reviewed by V. P. Snaith 

MR1756354 (2001d:19001)
Nakamura, Jinya(J-TOKYOGM)
On the Milnor $K$-groups of complete discrete valuation
fields. (English. English summary) 
Doc. Math. 5 (2000), 151--200 (electronic).
19D45 (11S70)
                                                                                
Summary: "For a discrete valuation field $K$, the unit group
$K^\times$ of $K$ has a natural decreasing filtration with respect to
the valuation, and the graded quotients of this filtration are given
in terms of the residue field. The Milnor $K$-group $K^M_q(K)$ is a
generalization of the unit group, and it also has a natural decreasing
filtration. However, if $K$ is of mixed characteristics and has an
absolute ramification index greater than one, then the graded
quotients of this filtration are not yet known except in some special
cases. The aim of this paper is to determine them when $K$ is an
absolutely tamely ramified discrete valuation field of mixed
characteristics $(0,p>2)$ with possibly imperfect residue
field. Furthermore, we determine the kernel of the Kurihara
$K^M_q$-exponential homomorphism from the differential module to the
Milnor $K$-group for such a field."

MR1948685 (2004c:11116)
Fukaya, Takako(J-TOKYOGM)
The theory of Coleman power series for $K\sb 2$. (English. English summary)
J. Algebraic Geom. 12 (2003), no. 1, 1--80.
11G55 (11F11 11S31)

This paper is a contribution to higher-dimensional local class field
theory following Kato's work. The author builds a "Coleman power
series" for the $K_2$ group of a local field of mixed characteristic
$(0,p)$ with imperfect residue field $k$ such that $[k \colon
k^p]=p$. The Coleman power series is an explicit reciprocity law for
any Lubin-Tate formal group in the case of a perfect residue
field. The case studied by the author corresponds to the
multiplicative group but in the case of imperfect residue field. In
the main part of the work, the author uses syntomic cohomology as
initiated by Kato to get his result.

MR1979649 (2004i:11057)
Kimura, Kenichiro(J-TSUKS-IM)
Elliptic units in $K\sb 2$. (English. English summary)
J. Number Theory 101 (2003), no. 1, 1--12.
11G16 (11G55)

In the paper the author constructs norm-compatible systems in $K_2$ of
function fields of some CM elliptic curves. In a more general setting
these systems were considered by A. J. Scholl \ref[in Galois
representations in arithmetic algebraic geometry (Durham, 1996),
379--460, Cambridge Univ. Press, Cambridge, 1998; MR1696501
(2000g:11057)]. In Section 2 of the paper the author applies the
theory of Coleman power series for $K_2$ of two-dimensional local
fields introduced by T. Fukaya \ref[J. Algebraic Geom. 12 (2003),
no. 1, 1--80; MR1948685 (2004c:11116)]. It is shown that this
construction leads to some elements of $K_2$ of certain power series
rings, which can be described explicitly. In Section 3 the author
shows that after some modification the systems are in the kernel of
tame symbol maps. He also describes the sufficient condition for the
images of the systems under the Chern class maps to be nonzero. In the
appendix, following Scholl's argument, the author shows that the
considered systems form Euler systems.

Reviewed by Piotr Kraso\'n 

MR0908217 (90d:28022)
Wagoner, J. B.(1-CA)
Markov partitions and $K\sb 2$.
Inst. Hautes Études Sci. Publ. Math. No. 65, (1987), 91--129.
28D05 (18F25 19C99 20B27 46L80 54H20 60J10)

Let $\scr S$ be a finite or countably infinite state space, let
$A\colon\scr S\times \scr S\to\{0,1\}$ be a zero-one matrix, and let
$\sigma_A\colon X_A\to X_A$ be the (two-sided) Markov shift associated
to $A$, i.e. $X_A=\{x=\{x_k\}\in\scr S^Z\colon\break A(x_k,x_{k+1})=1$
for all $k\in\bold Z\}$, which has a metric space structure, and
$(\sigma_A(x))_k=x_{k+1}$. Let ${\rm Aut}(\sigma_A)$ denote the
discrete group of all uniformly continuous homeomorphisms
$\alpha\colon X_A\to X_A$ which commute with $\sigma_A$ and such that
$\alpha^{-1}$ is uniformly continuous also. The complete structure of
${\rm Aut}(\sigma_A)$ is still unknown, although many partial results
have been obtained by \n G. Hedlund\en \ref[Math. Systems Theory 3
(1969), 320--375; MR0259881 (41 \#4510)] and more recently by \n
M. Boyle\en and \n W. Krieger\en \ref[Trans. Amer. Math. Soc. 302
(1987), no. 1, 125--149; MR0887501 (88g:54065)] among others.

The main aim of the interesting and innovative paper under review is
to obtain information about ${\rm Aut}(\sigma_A)$ by constructing
homomorphisms from it into other groups. Towards this end the author
gives $\scr P_A$, the set of topological Markov partitions for
$\sigma_A$ on $X_A$ (see \n W. Parry\en and \n S. Tuncel's\en
Classification problems in ergodic theory \ref[Cambridge Univ. Press,
Cambridge, 1982; MR0604728 (82g:28024)] for definitions), the
structure of a simplicial complex on which ${\rm Aut}(\sigma_A)$ acts
in a properly discontinuous fashion, and proves that $\scr P_A$ is
contractible and is locally compact when $A$ is finite. In doing so,
the author answers a question of \n M. Gromov\en concerning
combinatorial properties of $\scr P_A$. Using these results the author
first constructs a homomorphism from ${\rm Aut}(\sigma_A)$ to the
fundamental group $\pi_1(S(\overline{\scr E}),A)$, where
$S(\overline{\scr E})$ denotes a simplicial complex formed from the
set of shift equivalences of a certain class of nonnegative matrices
\ref[see R. Williams, Ann. of Math (2) 98 (1973), 120--153; errata;
MR0331436 (48 \#9769)]. In the case where $A$ is finite, using the
method of "inverting functors" the author then constructs an
isomorphism between $\pi_1(S(\overline{\scr E}),A)$ and ${\rm
Aut}(G(A),G(A)_+,S_A)$, the group of order preserving automorphisms of
the dimension group $G(A)$ associated to $A$ \ref[Krieger,
Invent. Math. 56 (1980), no. 3, 239--250; MR0561973
(81m:28018)]. Again in the case where $A$ is finite, generalizations
of the two maps above can be composed with a determinant map to
construct a commutative diagram $$\matrix {\rm Aut}(\sigma_A)\buildrel
K_A\over\longrightarrow K_2(F(t))\cr
\scriptstyle{\kappa_A}\textstyle\searrow\qquad\textstyle
\downarrow{\ssize\delta}\cr \qquad\bigoplus_{\scr P}F^*_{\scr
P}\endmatrix$$ where $F(t)$ is the field of rational functions over a
field $F, \scr P$ is a prime ideal of $F[t,t^{-1}]$, $F_{\scr
P}=F[t,t^{-1}]/\scr P,K_2$ is the algebraic $K$-group and $\delta$ is
the tame symbol of \n J. W. Milnor\en \ref[ Introduction to algebraic
$K$-theory, Ann. of Math. Stud., 72, Princeton Univ. Press, Princeton,
NJ, 1971; MR0349811 (50 \#2304)], and the image of $\kappa_A$ is
contained in the sum of those $F^*_{\scr P}$ where $\scr P$ divides
$\det(I-tA)$.

Other applications of the author's results are also given, to the
theory of random walks on a countably infinite discrete group, and to
Williams' theory of strong shift equivalence and Markov chains,
generalizing results of Williams \ref[op. cit.] and Parry and Tuncel
\ref[op. cit.]. The author has subsequently extended some of his
methods and has discussed in further detail their connection with the
Cuntz-Krieger $C^*$-algebras \ref[Krieger and J. Cuntz,
Invent. Math. 56 (1980), no. 3, 251--268; MR0561974 (82f:46073a)] in
another paper \ref[the author, Adv. in Math. 71 (1988), no. 2,
133--185].

Reviewed by Judith A. Packer

MR1220423 (94d:19003)
Barge, Jean(F-POLY)
Cocycle d'Euler et $K\sb 2$. (French. English summary) [Euler cocycle
and $K\sb 2$] 
$K$-Theory 7 (1993), no. 1, 9--16.
19C30 (11E81)

For every field $k$, a canonical 2-cocycle is constructed on the group
${\rm PSL}(2,k)$ with values in the Witt group $W(k)$. This yields a
natural homomorphism from $H_2({\rm SL}(2,k);Z)$ into the square
$I^2(k)$ of the fundamental ideal of $W(k)$. The main result is that
this homomorphism is a lift of the Milnor symbol which associates to a
symbol $\{x,y\}\in K_2(k)$ the class of the Pfister form
$\langle\langle x,y\rangle\rangle$ in $I^2(k)/I^3(k)$.

Reviewed by Jurgen Hurrelbrink

MR2114169 (2005k:19007)
Akhtar, Reza(1-MMOH)
Milnor $K$-theory of smooth varieties. (English. English summary)
$K$-Theory 32 (2004), no. 3, 269--291.
19D45 (14C25 14C35)

Let $X$ be a smooth scheme of finite type over a field $k$. Then one
has an isomorphism between Bloch's higher Chow groups and motivic
cohomology groups with coefficients in some abelian group $A$:
$$H_{\scr M}^{2q}(X,A(q))\simeq {\rm CH}^q(X,2q-p;A).$$ In the case
$X={\rm Spec}\,F, F$ a field and $A=\Bbb Z$, then $K_n^M(F)\simeq {\rm
CH}^n(F,n)$, where $K_*^M$ is Milnor's $K$-theory. Generalizing
results by M. Somekawa \ref[$K$-Theory 4 (1990), no. 2, 105--119;
MR1081654 (91k:11052)] and by W. Raskind and M. Spieß \ref[Compositio
Math. 121 (2000), no. 1, 1--33; MR1753108 (2002b:14007)], the author
extends the setting of Milnor $K$-theory from the case of fields to
that of a smooth projective variety in the following way: given a
geometrically integral quasi-projective variety $X$ which is smooth of
dimension $d$ over $k$, one defines "mixed $K$-groups" $K(k,\scr
{CH}_0(X),\Bbb G_m)$ (where $\Bbb G_m$ is the multiplicative group
scheme) such that for $X$ projective there is an isomorphism
$$K_s(k,\scr {CH}_0(X),\Bbb G_m)\simeq {\rm CH}^{d+s}(X,s).$$
Furthermore, for any smooth quasiprojective variety there is an
isomorphism $$H_{\rm {Zar}}^d(X,\scr K_{d+s}^M)\simeq {\rm
CH}^{d+s}(X,s),$$ where $\scr K_{d+s}^M$ is the sheafified version of
the Milnor $K$-groups.

The definition of the mixed $K$-groups $$K_s(k,\scr
{CH}_0(X_1),\dots,\scr {CH}_0(X_r);G_1,\dots,G_s)$$ for smooth
quasiprojective varieties $X_1,\dots ,X_r$ and semi-abelian varieties
$G_1,\dots ,G_s$ over the field $k$ is given in terms of the group
$$T=\bigoplus_{E/k\,{\rm finite}}{\rm CH}_0(X_1)_E\otimes \cdots
\otimes {\rm CH}_0(X_r)_E\otimes G_1(E)\otimes \cdots \otimes G_s(E)$$
and a subgroup $R$ of relations in $T$ which are similar to the
projection formula and the reciprocity law. Then one takes $K_s(k,\scr
{CH}_0(X),\Bbb G_m)$ to be the group $K_s(k,\scr
{CH}_0(X_1),\dots,\scr C\scr H_0(X_r);G_1,\dots,G_s)$ in the case
$r=1$ and $G_1=\dots=G_s=\Bbb G_m$, and $K_s(k,\Bbb G_m)$ in the case
$r=0$. The following isomorphism holds: $K_s(k;\Bbb G_m)\simeq
K_s^M(k)$.

Reviewed by Claudio Pedrini

MR1220423 (94d:19003)
Barge, Jean(F-POLY)
Cocycle d'Euler et $K\sb 2$. (French. English summary) [Euler cocycle
and $K\sb 2$] 
$K$-Theory 7 (1993), no. 1, 9--16.
19C30 (11E81)

For every field $k$, a canonical 2-cocycle is constructed on the group
${\rm PSL}(2,k)$ with values in the Witt group $W(k)$. This yields a
natural homomorphism from $H_2({\rm SL}(2,k);Z)$ into the square
$I^2(k)$ of the fundamental ideal of $W(k)$. The main result is that
this homomorphism is a lift of the Milnor symbol which associates to a
symbol $\{x,y\}\in K_2(k)$ the class of the Pfister form
$\langle\langle x,y\rangle\rangle$ in $I^2(k)/I^3(k)$.

Reviewed by Jurgen Hurrelbrink

K_2 of elliptic surface minus singular fibers and q-expansion of
Beilinson's regulator 
Authors: Masanori Asakura
Comments: 49 pages
Subj-class: Algebraic Geometry; Number Theory
MSC-class: 14J27; 19D45; 19F27

    The main result is to give a upper bound of the rank of dlog image
    of K_2 of elliptic surface minus singular fibers. This is based on
    a detailed study of K_2 of Tate curves over Z_p((q)). In
    particular Beilinson's regulator plays an important role. As a
    corollary, we have indecomposable parts of K_1 of elliptic surface
    with higher rank.

\
Paper: math.KT/0603241
Date: Fri, 10 Mar 2006 06:34:33 GMT   (20kb)

Title: Motivic interpretation of Milnor $K$-groups attached to Jacobian
 varieties
Authors: Satoshi Mochizuki
Comments: 29 pages
Subj-class: K-Theory and Homology; Algebraic Geometry
MSC-class: 19F15;11R58;14L10;14C25;14F42
\\
 In the paper M. Somekawa, {\it{On Milnor $K$-groups attached at semi-Abelian
varieties}}, K-theory, \textbf{4} (1990) p.105, Somekawa conjectures that his
Milnor K-group $K(k,G_1,...,G_r)$ attached to semi-abelian varieties
$G_1$,...,$G_r$ over a field $k$ is isomorphic to ${\rm
Ext}_{\mathcal{M}_k}^r(\mathbb{Z},G_1[-1] \otimes ... \otimes G_r[-1])$ where
$\mathcal{M}_k$ is a certain category of motives over $k$. The purpose of this
note is to give remarks on this conjecture, when we take $\mathcal{M}_k$ as
Voevodsky's category of motives ${\rm DM}^{eff}_{-}(k)$ .
\\ ( http://arXiv.org/abs/math/0603241 ,  20kb)

MR0743941 (85g:11110)
Bak, Anthony(D-BLF); Rehmann, Ulf(D-BLF)
$K\sb{2}$-analogs of Hasse's norm theorems.
Comment. Math. Helv. 59 (1984), no. 1, 1--11.
11R70 (16A54 18F25)

From the introduction: "The classical norm theorem of Hasse for fields
and the classical norm theorem of Hasse and Schilling for simple
algebras can be stated in the language of algebraic $K$-theory in
terms of the functor $K\sb 1$. The goal of this paper is to show that
$K\sb 2$-analogues of these results are valid. It turns out that
whereas the classical norm theorem for fields is valid essentially
only for cyclic extensions of global fields, its $K\sb 2$-analogue is
valid for all finite extensions of global fields."

Main results: Theorem 2: If $L$ is a finite extension of a global
field $K$, and $\Sigma\sb {L/K}$ is the set of all real places of $K$
such that all their extensions to $L$ are complex, then the following
sequence is exact: $$K\sb 2L\,{\buildrel {N\sb {L/K}}\over{\hbox to
20pt{\rightarrowfill}}}\,K\sb 2K\,{\buildrel\lambda\over{\hbox to
20pt{\rightarrowfill}}}\coprod\sb {v\in\Sigma\sb {L/K}}µ(K\sb v)\to
1,$$ where $N\sb {L/K}$ is the transfer mapping in $K$-theory,
$\lambda=\coprod\sb {v\in\Sigma\sb {L/K}}\lambda\sb v$, and
$\lambda\sb v$ is the norm residue symbol. Theorem 3: If $K$ is a
global field, $D$ is a finite $K$-central division algebra, and
$\Sigma\sb {D/K}$ is the set of all real places of $K$ not splitting
$D$, then the following sequence is exact: $$K\sb 2D\,{\buildrel {N\sb
{D/K}}\over{\hbox to 20pt{\rightarrowfill}}}\,K\sb
2K\,{\buildrel\lambda\over{\hbox to 20pt{\rightarrowfill}}}\coprod\sb
{v\in\Sigma\sb {D/K}}µ(K\sb v),$$ where $N\sb {D/K}$ is the reduced
norm mapping in $K$-theory and $\lambda$ is as above.

The authors reduce the proof of Theorem 3 to Theorem 2, and they prove
Theorem 2 by number-theoretic considerations. They do not use Tate's
results on relations between $K$-theory and Galois cohomology. The
appendix contains a straightforward proof of Theorem 3 given by
O. Gabber.

Reviewed by J. Browkin

MR2143560 (2006b:19004)
Kurihara, Masato(J-TOKYM)
On the structure of Milnor $K$-groups of certain complete discrete
valuation fields.  
J. Théor. Nombres Bordeaux 16 (2004), no. 2, 377--401.
19D45 (11S70)
 
Let $K$ be a discrete valuation field with ring of integers $O_K$,
maximal ideal $m_K$ and residue field $F$. The higher Milnor
$K$-groups $K_q(K)$ ($q>0$) have a natural filtration by subgroups
$U^iK_q(K)$, which are generated by $\{1+m_K^i,K^*,\dots,K^*\}$. The
graded quotients $$ gr^iK_q(K)=U^iK_q(K)/U^{i+1}K_q(K) $$ are known in
the equal characteristic case \ref[S. Bloch, Inst. Hautes Études
Sci. Publ. Math. No. 47 (1977), 187--268 (1978); MR0488288
(81j:14011)] and in the mixed characteristic case in a certain range
\ref[S. J. Bloch and K. Kat\=o, Inst. Hautes Études
Sci. Publ. Math. No. 63 (1986), 107--152; MR0849653 (87k:14018)].

The author determines the graded quotients $gr^iK_2(K)$ for Milnor's
$K_2$ for the following typical example in mixed characteristic, which
is a discrete valuation field of type II in the terminology of his
paper "On two types of complete discrete valuation fields"
\ref[Compositio Math. 63 (1987), no. 2, 237--257; MR0906373
(88k:11090)]. Let $K_0$ denote the field of fractions of the
completion of the localization of $\Bbb Z_p[T]$ at the prime ideal
$(p)$, and let $K=K_0(\root p\of {pT})$.

In the interesting range $i>p+1$ not covered by the results of
Bloch-Kato the author shows that $$ gr^iK_2(K)\cong\cases 0&\quad{\rm
if} i {\rm is prime to} p,\\ F/F^p&\quad{\rm if} i=2p,\\
F^{p^{n-2}}&\quad{\rm if} i=np, n\geq3.\endcases $$

As a consequence, the author obtains the interesting result that $K$
does not have a cyclic, totally ramified extension of degree $p^3$.

Reviewed by Manfred Kolster 

MR0906373 (88k:11090)
Kurihara, Masato(J-TOKYOS)
On two types of complete discrete valuation fields.
Compositio Math. 63 (1987), no. 2, 237--257.
11S70 (18F25 19D45 19F05)

Complete discrete valuation fields of mixed characteristics $(0,p)$
are considered. They are divided into two types. There is much
difference between the structures of Milnor $K$-groups and between the
abelian extensions for these two types.

Reviewed by L. N. Vaserstein

\bye